\documentclass[a4paper]{amsart}

\usepackage{times}
\usepackage[vlined,algonl,ruled,algosection,noresetcount]{algorithm2e}
\usepackage{epsfig}
\usepackage[T1]{fontenc}
\usepackage[latin1]{inputenc}
\usepackage[arrow,curve,matrix]{xy}

\newcommand{\PreprintAndPublication}[2]{#1} 

\DeclareFontFamily{OMS}{rsfs}{\skewchar\font'60}
\DeclareFontShape{OMS}{rsfs}{m}{n}{<-5>rsfs5 <5-7>rsfs7 <7->rsfs10 }{}
\DeclareSymbolFont{rsfs}{OMS}{rsfs}{m}{n}
\DeclareSymbolFontAlphabet{\scr}{rsfs}

\renewcommand{\P}{\mathbb{P}}
\newcommand{\C}{\mathbb{C}}
\renewcommand{\O}{\mathcal{O}}

\newcommand\resto[1]{\hbox{\hbox{${\vert{}_{{#1}}}$}}}

\newcommand{\sA}{\scr{A}} 
\newcommand{\sL}{\scr{L}} 
\newcommand{\sR}{\scr{R}} 

\newcommand{\bN}{\mathbb{N}}

\newcommand{\bP}{\mathbb{P}}
\newcommand{\bQ}{\mathbb{Q}}

\DeclareMathOperator{\bir}{bir}
\DeclareMathOperator{\Chow}{Chow}
\DeclareMathOperator{\codim}{codim}
\DeclareMathOperator{\coker}{coker}
\DeclareMathOperator{\discrep}{discrep}
\DeclareMathOperator{\Hom}{Hom}

\DeclareMathOperator{\Pic}{Pic}
\DeclareMathOperator{\red}{red}
\DeclareMathOperator{\Sing}{Sing}
\DeclareMathOperator{\Sym}{Sym}
\DeclareMathOperator{\supp}{supp}

\DeclareMathOperator{\Var}{Var}
\DeclareMathOperator{\Spec}{Spec}

\DeclareMathOperator{\id}{id}

\newcounter{thisthm}
\newcommand{\ilabel}[1]{\newcounter{#1}\setcounter{thisthm}{\value{thm}}\setcounter{#1}{\value{enumi}}}
\newcommand{\iref}[1]{(\thesection.\the\value{thisthm}.\the\value{#1})}

\newtheoremstyle{bozont}{3pt}{3pt}%
     {\itshape}
     {}
     {\bfseries}
     {\ }
     {.5em}
     {\thmname{#1}\thmnumber{ #2.}\thmnote{{ \rm #3.}}}
\theoremstyle{bozont}
\newtheorem{thm}{Theorem}[section]
\newtheorem{conj}[thm]{Conjecture}
\newtheorem{defn}[thm]{Definition}

\newtheorem{assumption}[thm]{Assumption}
\numberwithin{equation}{thm}
\numberwithin{figure}{section}
\theoremstyle{plain}    
\newtheorem{cor}[thm]{Corollary}
\newtheorem{lem}[thm]{Lemma}

\newtheorem{num}[thm]{}
\theoremstyle{plain}    
\newtheorem{prop}[thm]{Proposition}
\newtheorem{proclaim-special}[thm]{\specialthmname}

\theoremstyle{remark}
\newtheorem{fact}[thm]{Fact}
\newtheorem{rem}[thm]{Remark}
\newtheorem{obs}[thm]{Observation}
\newtheorem{subrem}[equation]{Remark}
\newtheorem{subclaim}[equation]{Claim} 
\newtheorem*{claim*}{Claim} 
\newtheorem{notation}[thm]{Notation}
\newtheorem{construction}[thm]{Construction}

\newtheoremstyle{bozont-remark}{3pt}{3pt}%
     {}
     {}
     {\it}
     {.}
     {.5em}
     {\thmname{#1}\thmnumber{ #2}: \thmnote{\sc #3}}
\theoremstyle{bozont-remark}

\def\factor#1.#2.{\left. \raise 2pt\hbox{$#1$} \right/\hskip -2pt\raise -2pt\hbox{$#2$}} 

\renewcommand\thesubsection{\thesection.\Alph{subsection}}
\makeatletter
\renewcommand\thesubsubsection{\thesubsection.\@arabic\c@subsubsection}
\makeatother

\newenvironment{narrow}{
  \begin{list}{}%
    {
        \setlength\leftmargin{.3\parindent}
    }
  \item[]\vskip-12pt}{\end{list}}
\newenvironment{nequation}{\begin{minipage}{\textwidth}\begin{equation}}
    {\end{equation}\end{minipage}\medskip}
\newenvironment{newequation}{\numberwithin{equation}{section}
  \setcounter{equation}{\value{thm}} %
    \begin{equation}}
    {\end{equation}\addtocounter{thm}1\setcounter{equation}0
    \numberwithin{equation}{thm}} 

\newenvironment{enumerate-c}{
  \begin{enumerate}
    \setcounter{enumi}{\value{equation}}}
  {\setcounter{equation}{\value{enumi}}\end{enumerate}}

\begin{document}

\title{Families of canonically polarized varieties over surfaces}

\date{\today}

\author{Stefan Kebekus}
\author{S\'andor J.\ Kov\'acs}

\thanks{Both authors were supported in part by the priority program
  ``Globale Methoden in der komplexen Geometrie'' of the Deutsche
  Forschungsgemeinschaft, DFG.  S\'andor Kov\'acs was supported in
  part by NSF Grant DMS-0092165 and a Sloan Research Fellowship. A
  part of this paper was worked out while Stefan Kebekus visited the
  Korea Institute for Advanced Study.  He would like to thank Jun-Muk
  Hwang for the invitation.}

\address{Stefan Kebekus, Mathematisches Institut, Universit\"at zu
  K\"oln, Weyertal 86--90, 50931 K\"oln, Germany}
\email{\href{mailto:stefan.kebekus@math.uni-koeln.de}{stefan.kebekus@math.uni-koeln.de}}
\urladdr{\href{http://www.mi.uni-koeln.de/~kebekus}{http://www.mi.uni-koeln.de/$\sim$kebekus}}

\address{\noindent S\'andor Kov\'acs, University of Washington, Department
  of Mathematics, Box 354350, Seattle, WA 98195, U.S.A.}
\email{\href{mailto:kovacs@math.washington.edu}{kovacs@math.washington.edu}}
\urladdr{\href{http://www.math.washington.edu/~kovacs}{http://www.math.washington.edu/$\sim$kovacs}}

\date{\today}

\begin{abstract}
  Shafarevich's hyperbolicity conjecture asserts that a family of
  curves over a quasi-projective 1-dimensional base is isotrivial
  unless the logarithmic Kodaira dimension of the base is positive.
  More generally it has been conjectured by Viehweg that the base of a
  smooth family of canonically polarized varieties is of log general
  type if the family is of maximal variation.  In this paper, we
  relate the variation of a family to the logarithmic Kodaira
  dimension of the base and give an affirmative answer to Viehweg's
  conjecture for families parametrized by surfaces.
\end{abstract}

\maketitle

\setcounter{tocdepth}{1}
\tableofcontents

\section{Introduction}

Let $B^\circ$ be a smooth quasi-projective complex curve and $q>1$ a positive
integer.  Shafarevich conjectured \cite{Shaf63} that the set of
non-iso\-trivial families of smooth projective curves of genus $q$
over $B^\circ$ is finite. Shafarevich further conjectured that if the
logarithmic Kodaira dimension, for a definition see below, satisfies
$\kappa(B^\circ) \leq 0$, then no such families exist.  This
conjecture, which later played an important role in Faltings' proof of
the Mordell conjecture, was confirmed by Parshin \cite{Parshin68} for
$B^\circ$ projective and by Arakelov \cite{Arakelov71} in general. We
refer the reader to the survey articles \cite{Viehweg01} and
\cite{Kovacs03c} for a historical overview and references to related
results.

It is a natural and important question whether similar statements hold for families
of higher dimensional varieties over higher dimensional bases.  Families over a curve
have been studied by several authors in recent years and they are now fairly well
understood---the strongest results known were obtained in \cite{Vie-Zuo01,VZ02}, and
\cite{Kovacs02}.  For higher dimensional bases, however, a complete picture is still
missing and subvarieties of the corresponding moduli stacks are not well understood.
As a first step toward a better understanding, Viehweg proposed the following:

\begin{conj}[\protect{\cite[6.3]{Viehweg01}}]\label{conj:viehweg}
  Let $f^\circ: X^\circ \to S^\circ$ be a smooth family of canonically
  pol\-arized varieties. If $f^\circ$ is of maximal variation, then
  $S^\circ$ is of log general type.
\end{conj}

We briefly recall the relevant definitions, as they will also be
important in the statement of our main result. The first is the
variation, which measures the birational non-isotriviality of a
family.

\begin{defn}\label{def:var}
  Let $f: X \to S$ be a projective family over an irreducible base $S$ 
defined over an algebraically closed field $k$
  and let $\overline{k(S)}$ denote the algebraic closure of the function
  field of $S$.  The variation of $f$, denoted by $\Var f$, is defined
  as the smallest integer $\nu$ for which there exists a subfield $K$
  of $\overline{k(S)}$, finitely generated of transcendence degree
  $\nu$ over $k$ and a $K$-variety $F$ such that
  $X\times_S\Spec\overline{k(S)}$ is birationally equivalent to
  $F\times_{\Spec K}\Spec\overline{k(S)}$.
\end{defn}

\begin{subrem}
  In the setup of Definition~\ref{def:var}, if the fibers are canonically polarized
  complex varieties, moduli schemes are known to exist, and the variation is the same
  as either the dimension of the image of $S$ in moduli, or the rank of the
  Kodaira-Spencer map at the general point of $S$.
\end{subrem}

\begin{defn}
  Let $S^\circ$ be a smooth quasi-projective variety and $S$ a smooth
  projective compactification of $S^\circ$ such that $D := S \setminus
  S^\circ$ is a divisor with simple normal crossings. The logarithmic
  Kodaira dimension of $S^\circ$, denoted by $\kappa(S^\circ)$, is
  defined to be the Kodaira-Iitaka dimension, $\kappa(S,D)$, of the
  line bundle $\O_S(K_S + D) \in \Pic(S)$.
  The variety $S^\circ$ is called of \emph{log general type} if
  $\kappa(S^\circ)=\dim S^\circ$, i.e., the divisor $K_S+D$ is big.
\end{defn}

\begin{subrem}
  It is a standard fact in logarithmic geometry that a
  compactification $S$ with the described properties exists, and that the logarithmic
  Kodaira dimension $\kappa(S^\circ)$ does not depend on the choice of the
  compactification $S$.
\end{subrem}

\subsection{Statement of the main result}

Our main result describes families of canonically polarized varieties over
quasi-projective surfaces. We relate the variation of the family to the logarithmic
Kodaira dimension of the base and give an affirmative answer to Viehweg's Conjecture
\ref{conj:viehweg} for families over surfaces.

\begin{thm}\label{thm:main}
  Let $S^\circ$ be a smooth quasi-projective complex surface and $f^\circ:
  X^\circ \to S^\circ$ a smooth non-isotrivial family of canonically
  polarized complex varieties. Then the following holds.
  \begin{enumerate}     
  \item If $\kappa(S^\circ)=-\infty$, then $\Var(f^\circ) \leq 1$.
    
  \item If $\kappa(S^\circ) \geq 0$, then $\Var(f^\circ) \leq
    \kappa({S^\circ})$.
  \end{enumerate}
  In particular, Viehweg's Conjecture holds for families over surfaces,
\end{thm}

\noindent
For the special case of $\kappa(S^\circ)=0$, this statement was proved by Kov\'acs
\cite[0.1]{Kovacs97a} when $S^\circ$ is an abelian variety and more generally by Viehweg
and Zuo \cite[5.2]{VZ02} when $T_S(-\log D)$ is weakly positive.

A slightly weaker statement holds for families of minimal varieties,
see Section~\ref{sec:complements} below.  In a forthcoming paper we
will give a more precise geometric description of $f^\circ$ in the
case of $\kappa(S^\circ) \leq 1$.

\begin{rem}
  Notice that in the case of $\kappa(S^\circ)=-\infty$ one cannot
  expect a stronger statement. For an easy example take any
  non-isotrivial smooth family of canonically polarized varieties over
  a curve $g: Z\to C$, set $X : = Z\times \bP^1$, $S^\circ := C\times
  \bP^1$, and let $f^\circ := g \times \id_{\bP^1}$ be the obvious
  morphism. Then we clearly have $\kappa(S^\circ)=-\infty$ and
  $\Var(f)=1$.
\end{rem}

\noindent
In view of Theorem~\ref{thm:main}, we propose the following generalization of
Viehweg's conjecture.

\begin{conj}
  Let $f^\circ: X^\circ \to S^\circ$ be a smooth family of canonically polarized
  varieties. Then either $\kappa(S^\circ)=-\infty$ and $\Var(f^\circ)<\dim S^\circ$,
  or $\Var(f^\circ)\leq \kappa(S^\circ)$.
\end{conj}

\subsection{Outline of the paper}

Throughout the paper we work over $\mathbb C$, the field of complex numbers.

The paper is divided into two parts. In the first part comprising
Sections~\ref{sec:loggeometry} and \ref{subsec:rel-min} we recall and establish
techniques that might be of independent interest.  Section~\ref{sec:loggeometry}
summarizes results in logarithmic geometry and logarithmic deformation theory. In
Section~\ref{subsec:rel-min} we consider logarithmic pairs $(S,D)$ where $S$ is a
birationally ruled surface, and construct a sequence of blowing down $(-1)$-curves
that can be used to simplify the self-intersection graph of the boundary $D$.

In the second part of the paper we employ these techniques in order to
prove Theorem~\ref{thm:main}. After the notation is set up in
Section~\ref{sec:setup}, we consider the cases where the logarithmic
Kodaira dimension of $S^\circ$ is 1, 0 or $-\infty$ in
Sections~\ref{sec:k1}--\ref{sec:Kinfinity}, respectively.

\subsection{Acknowledgements}

The authors would like to thank J\'anos Koll\'ar for calling their
attention to an error in a previous version of the paper and for
immediately suggesting a correction. The authors would also like to
thank the referee, who suggested to mention the generalizations in
Section~\ref{sec:complements}.

\part{TECHNIQUES}

\section{Logarithmic geometry}
\label{sec:loggeometry}

Throughout the current section, let $S$ be a smooth projective variety and $D \subset
S$ a reduced divisor with simple normal crossings. As follows, we recall a number of
facts concerning this setup and include proofs wherever we could not find an adequate
reference.

\subsection{The Logarithmic Minimal Model Program}
\label{sec:logMMP}

If $S$ is a surface and the logarithmic Kodaira dimension
$\kappa(K_S+D)$ is non-negative, we will frequently need to consider
the $(S,D)$-logarithmic minimal model program, which is briefly
recalled here. The reader is referred to \cite{KM98} for the relevant
definitions, for proofs and for a full discussion.

\begin{fact}[Logarithmic Minimal Model Program,
  \protect{\rm \cite[(3.47)]{KM98}}]\label{fact:lmmp} If $\dim S=2$ and
  $\kappa(K_S+D) \geq 0$, there exists a birational morphism $\phi : S
  \to S_{\lambda}$ from $S$ to a normal surface $S_{\lambda}$ such
  that
  \begin{enumerate}
  \item The morphism $\phi$ is the composition of finitely many
    log contractions.
  
  \item If we set $D_{\lambda} := \phi(D)$ to be the cycle-theoretic image divisor,
    then
    \begin{enumerate}
    \item The pair $(S_{\lambda}, D_{\lambda})$ has only \emph{dlt} singularities and
      $S_{\lambda}$ itself is $\mathbb Q$-factorial \cite[(3.36), (3.44)]{KM98}. In
      particular, $S_{\lambda}$ has only quotient singularities.
      
    \item The log canonical divisor $K_{S_{\lambda}} + D_{\lambda}$ is nef.
    
    \item The log Kodaira dimension remains unchanged,
      $$
      \kappa(K_{S_{\lambda}} + D_{\lambda} ) = \kappa(K_S + D ).
      $$
    \end{enumerate}
  \end{enumerate}
\end{fact}

\PreprintAndPublication{
\begin{rem}\label{rem:singsinMMP}
  In the discussion of the minimal model program one needs to consider several
  classes of singularities. The large number of notions, and the fact that the
  definitions found in the literature are not always obviously equivalent makes the
  field somewhat difficult to navigate for the outsider. For the reader's
  convenience, we briefly indicate how that fact that $S_{\lambda}$ has only quotient
  singularities follows from the assumption that $(S_{\lambda}, D_{\lambda})$ has
  only $\bQ$-factorial \emph{dlt} singularities:
  
  By \cite[(2.42)]{KM98}, if $x \in S_\lambda$ is any point, then
  either $x$ is a smooth point of $S_\lambda$, or $(S_{\lambda},
  D_{\lambda})$ is \emph{plt} at $x$. We can thus assume without loss
  of generality that $(S_{\lambda}, D_{\lambda})$ is \emph{plt}
  everywhere, i.e.~that
  $$
  \discrep(S_\lambda, D_\lambda) > -1.
  $$
  By \cite[(2.27)]{KM98}, $\discrep(S_\lambda, 0) \geq \discrep(S_\lambda,
  D_\lambda) > -1$. By definition, \cite[(2.34)]{KM98}, this means that $S_\lambda$ is
  \emph{log terminal}. The classification of \emph{log terminal} surface
  singularities, \cite[(4.18)]{KM98} then gives the claim.
\end{rem}
}{}

\begin{rem}
  We remark that the support of $D_{\lambda}$ is generally \emph{not} equal to the
  image $\phi(D)$, as it may well happen that $\phi(D)$ contains isolated points
  which do not appear in the cycle-theoretic image.  This observation will later
  become important in Section~\ref{sec:reduction-to-birationally-ruled} and in the
  proof of Proposition~\ref{prop:uniLK0}.
\end{rem}

\begin{fact}[Logarithmic Abundance Theorem in Dimension 2,
  \protect{\cite[(3.3)]{KM98}}]\label{fact:abundance} 
  The linear system $|n(K_{S_{\lambda}}+D_{\lambda})|$ is basepoint-free for
  sufficiently large and divisible $n \in \mathbb N$.
\end{fact}

\subsection{Logarithmic deformation theory}

In Sections~\ref{sec:K0} and \ref{sec:Kinfinity} we will have to deal with families
of curves on $S$ that intersect the boundary divisor $D$ in one or two points. In
counting these points, intersection multiplicity does not play any role, but the
number of local analytic branches of the curves does.  More precisely, we use the
following definition.

\begin{defn}\label{defn:nptintersect}
  Let $X$ be an algebraic variety, $E \subset X$ an algebraic set, and $\ell \subset
  X$ a reduced proper curve with normalization $\nu : \widetilde \ell \to \ell$. We
  say that ``$\ell$ intersects $E$ in $d$ points'' if the preimage $\nu^{-1}(E)$ is
  supported on exactly $d$ closed points of $\widetilde \ell$.
\end{defn}

\PreprintAndPublication{
  \begin{figure}[htbp]
    \centering
    \begin{picture}(10.5, 2)
      \put(0.0, -.5){\includegraphics[height=2.5cm]{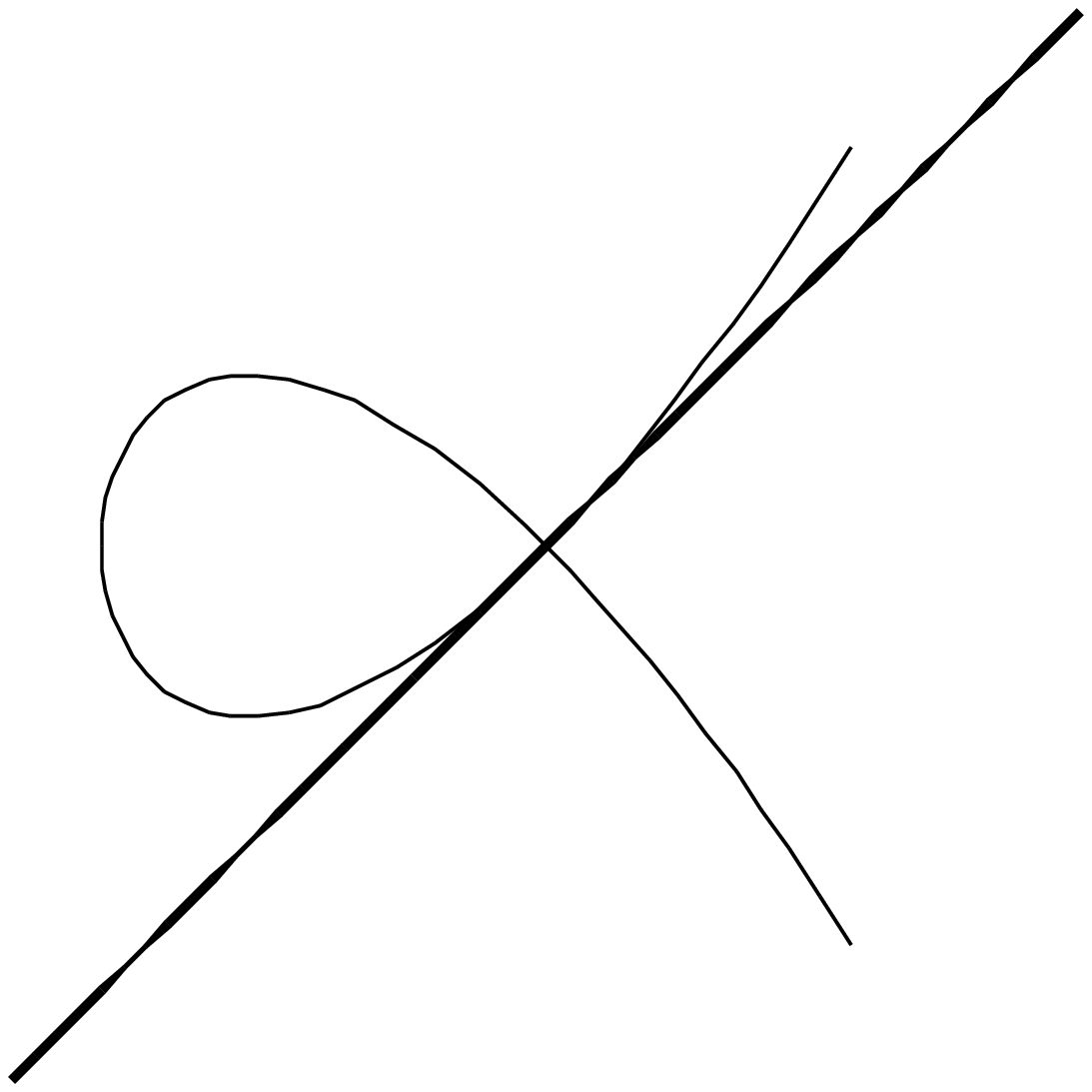}}
      \put(0.6, -.2){$E$}
      \put(2.0, 0.2){$\ell$}
      \put(3.5, 1.0){\parbox{7cm}{In the sense of Definition~\ref{defn:nptintersect},
          $\ell$ intersects $E$ in two points.}} 
    \end{picture}
    \caption{Number of intersection points}
  \end{figure}
}{}

\begin{rem}\label{rem:liftingOfCurves}
  Suppose we are given a proper birational morphism $\phi: X \to X'$, an algebraic
  set $E \subset X$ and a family of curves $\ell'_t \subset X'$ that intersect
  $\phi(E)$ in exactly $d$ points.  Assume further that none of the $\ell'_t$ is
  contained in the set of fundamental points of $\phi^{-1}$. Then the strict
  transforms give a (possibly disconnected) family $\ell_t$ of curves on $X$ that
  intersect $E$ in no more than $d$ points. If $E$ contains the $\phi$-exceptional
  locus, then the strict transforms intersect $E$ in exactly $d$ points.
\end{rem}

For our applications, we need to consider a family $\ell_t$ of rational curves in $S$
that intersect $D$ in two points.  Our aim in this section is to discuss an algebraic
parameter space for such curves. The construction is based on the observation that
for any such curve $\ell_t$ there exists a morphism $\nu_t : \P^1 \to \ell_t \subset
S$ such that $\nu_t^{-1}(D)$ is supported exactly on the points $x_0 := [0:1]$ and
$x_{\infty} := [1:0]$. Therefore, it makes sense to consider the space
$$
\mathcal H := \bigl\{ f \in \Hom_{\bir}(\P^1, S) \,|\, f^{-1}(D)
\text{ is supported exactly on $x_0$ and $x_{\infty}$} \bigr\}
$$
with the obvious structure as a closed, but possibly non-reduced
subscheme of $\Hom_{\bir}(\P^1, S)$, the space of generically
injective morphisms $\P^1 \to S$.

The space $\mathcal H$, and its infinitesimal structure has been studied in
\cite{KMcK} using a slightly different language. We recall some of their results here
and include proofs wherever we had difficulties to follow the original arguments.

\begin{prop}\label{prop:tgtspace}
  If $f \in \mathcal H$ is any closed point, then the Zariski tangent
  space to $\mathcal H$ at $f$ is canonically isomorphic to
  $$
  T_{\mathcal H}\resto f \simeq H^0\bigl(\P^1, f^*(T_S(-\log D))\bigr).
  $$
\end{prop}
\begin{proof}
  Let $\scr D=\{D_1, \dots, D_n\}$ be the set of irreducible components of $D$,
  $E_i:= f^*(D_i)$ the associated Cartier divisors on $\P^1$, and let $\scr E=\{E_1 ,
  \dots, E_n\}$.  Notice that all of the $E_i$ are supported on $x_0$ and $x_\infty$.
  If $\mathcal H_f \subset \mathcal H$ is the connected component that contains $f$,
  then $\mathfrak H = \Hom(\P^1, S, \scr E \subset \scr D)$, defined in
  \cite[Sect.~5]{KMcK}, is a subscheme of $\mathcal H$ that contains $\mathcal H_f$.
  Hence, the claim follows from \cite[(5.3)]{KMcK}.
\end{proof}

It is well known in the theory of rational curves on algebraic varieties that if $S$
is a uniruled manifold and $\ell$ is a rational curve that passes through a very
general point of $S$, then $\ell$ is free, and its deformations avoid any given
subset of codimension $\geq 2$.  More precisely, for any given subset $E \subset S$
with $\codim_S E \geq 2$ there exists a deformation $\ell'$ of $\ell$ that does not
intersect $E$. We show that a similar property holds for $\mathcal H$.

\begin{prop}[Small Set Avoidance]\label{prop:finitesetavoidance}
  Let $\mathcal H' \subset \mathcal H$ be an irreducible component
  such that the associated curves dominate $S$. If $M \subset
  S\setminus D$ is any closed set of $\codim_S M \geq 2$, then there
  exists a non-empty open set $\mathcal H'_0 \subset \mathcal H'$ such
  that for all $f' \in \mathcal H'_0$ the image does not intersect
  $M$, i.e.
  $$
  M \cap f'(\P^1) = \emptyset.
  $$
\end{prop}

The proof of Proposition~\ref{prop:finitesetavoidance}, which we give
on page~\pageref{proof:finitesetavoidance} at the end of this section
is based on a number of results that we prove first. We start with an
estimate for the dimension of $\mathcal H$ that we formulate and prove
in the next two lemmata.

\begin{defn}
  Let $\pi:X\to Y$ be a finite surjective morphism of degree $d$. The
  \emph{set-theoretic branch locus} of $\pi$ is the set of points in $Y$ whose
  set-theoretic preimage contains strictly less than $d$ points.
\end{defn}

\begin{lem}\label{lem:branching}
  Let $H$ be an irreducible variety and $D \subset \P^1 \times H$ an irreducible
  subvariety such that $\pi_2\resto D:D \to H$ is a finite surjective morphism of
  degree $d$ with set-theoretic branch locus $B$.  Then either $B = \emptyset$, or
  $B$ is a closed subvariety of pure codimension 1.
\end{lem}
\begin{proof}
  Performing a base change, if necessary, we can assume without loss of generality
  that $H$ is normal. The variety $D$ is then a well-defined family of algebraic
  cycles in the sense of \cite[I.3.10]{K96}, and therefore yields a morphism
  $$
  \gamma : H \to \Chow^d(\P^1) = \factor\P^1 \times \cdots \times
  \P^1.\text{permutation}. \simeq \P(\Sym^d {\mathbb A^2}).
  $$ 
  If $\Delta \subset \P(\Sym^d {\mathbb A^2})$ is the discriminant divisor,
  i.e. the branch locus of the morphism
  $$
  \P^1 \times \cdots \times \P^1 \longrightarrow \factor\P^1 \times
  \cdots \times \P^1.\text{permutation}. \simeq \P(\Sym^d {\mathbb A^2}),
  $$
  then the morphism $D \to H$ is branched at a point $\eta \in H$
  iff $\gamma(\eta) \in \Delta$. But since $\Delta \subset \P(\Sym^d
  {\mathbb A^2})$ is Cartier, the claim follows.
\end{proof}

The proof of Lemma~\ref{lem:branching} shows, after passing to the normalization,
that the branch locus $B$ is even a Cartier-divisor, but we will not need this
observation here.  The proposed estimate for the dimension of $\mathcal H$ then goes
as follows.

\begin{lem}[\protect{\cite[5.1, 5.3]{KMcK}}]\label{lem:dimofH}
  If $\eta \in \mathcal H$ is any point, then
  $$
  \dim_{\eta} \mathcal H \geq \dim_{\eta} \Hom_{\bir}(\P^1, S) -
  \underbrace{\deg_{\P^1} \eta^*(\O_S(D))}_{=:d}.
  $$
\end{lem}
\begin{proof}
  Let $H \subset \Hom_{\bir}(\P^1, S)$ be an irreducible component through $\eta$
  which is of maximal dimension. We will prove Lemma~\ref{lem:dimofH} by an inductive
  construction of a subvariety that contains $\eta$, satisfies the dimension bound,
  and is contained in $\mathcal H$. More precisely, we claim the following.

  \begin{narrow}
    \begin{subclaim}\label{dimension-estimate}
      There exists a sequence of subvarieties 
      $$
      H = H^{(0)} \supset H^{(1)} \supset \cdots \supset H^{(d-2)} \ni \eta
      $$
      such that $\codim_H H^{(i)} = i$, and such that for general
      closed points $f^{(i)} \in H^{(i)}$, we have $\#
      \left(f^{(i)}\right)^{-1}(D) \leq d -i$
    \end{subclaim}
    \begin{proof}[Proof] 
      We prove the claim inductively, using the index $i$ of $H^{(i)}$.  To start the
      induction, consider $i=0$. It is clear from intersection theory that if
      $f^{(0)} \in H$ is a general closed point, then $\# \bigl(f^{(0)}\bigr)^{-1}(D)
      \leq d$.  For the inductive step, assume that the subvariety $H^{(i)}$ is
      already constructed.  Consider the universal morphism $\mu_i : \P^1 \times
      H^{(i)} \to S $ and the reduced preimage
      $$
      D^{(i)} := \mu_i^{-1}(D)_{\red} \subset \P^1 \times H^{(i)}.
      $$
      By induction there exists an open set $\eta\in H^{(i)}_\circ\subseteq
      H^{(i)}$ such that $D^{(i)}_\circ=D^{(i)}\cap \pi_2^{-1}H^{(i)}_\circ$ surjects
      finitely onto $H^{(i)}_\circ$ with at most $d-i$ sheets.  Observe, that as long
      as $d-i>2$, $\eta$ will be in the set-theoretic branch locus of $\pi_2\resto
      {D^{(i)}_\circ} :D^{(i)}_\circ \to H^{(i)}_\circ$. If we set
      $$
      H^{(i+1)} := \text{closure of one component of the set-theoretic branch
        locus of } \pi_2\resto {D^{(i)}_\circ},
      $$
      then by Lemma~\ref{lem:branching} $\dim H^{(i+1)} = \dim H^{(i)}-1$, and a
      general point of $H^{(i+1)}$ has at most $d-(i+1)$ preimages on $D^{(i)}$.
      Claim~\ref{dimension-estimate} then follows.
    \end{proof}
  \end{narrow}
  
  \noindent
  According to Claim~\ref{dimension-estimate} there exists a subvariety containing
  $\eta$, $H^{(d-2)} \subseteq H$, of dimension $\dim H^{(d-2)} = \dim_{\eta} H - d +
  2$.  Since $\mu_{d-2}^{-1}(D)$ is a Cartier divisor on $\P^1 \times H^{(d-2)}$, the
  non-empty subvarieties
  \begin{align*}
    H^{(d-2)}_0 & := \pi_2 \left(\mu_{d-2}^{-1}(D)_{\red} \cap
      \{x_0\}\times H^{(d-2)} \right) \subseteq H^{(d-2)}, \text{ and } \\
    H^{(d-2)}_{0,\infty} & := \pi_2 \left(\mu_{d-2}^{-1}(D)_{\red} \cap
      \{x_\infty\}\times H^{(d-2)}_0 \right) \subseteq H^{(d-2)}_0
  \end{align*}
  each contain $\eta$ and have codimension at most 1 in one another. In other words,
  we have
  \begin{equation}
    \label{eq:dimh0inf}
    \dim H^{(d-2)}_{0,\infty} \geq \dim H^{(d-2)} -2 = \dim_{\eta} H - d. 
  \end{equation}
  It follows from Claim~\ref{dimension-estimate} that for all closed points $f \in
  H^{(d-2)}$, the associated morphism $f : \P^1 \to S$ satisfies $\#f^{-1}(D) \leq
  2$. Because $f$ is contained in $H^{(d-2)}_{0,\infty}\subseteq H^{(d-2)}_0$, we
  also have $f(x_0) \in D$ and $f(x_{\infty}) \in D$, respectively. In summary, we
  have seen that $H^{(d-2)}_{0,\infty} \subseteq \mathcal H$, which combined with
  \eqref{eq:dimh0inf} proves Lemma~\ref{lem:dimofH}.
\end{proof}

We note that a more detailed analysis of the construction could be
used to show that $\mathcal H$ is a local complete intersection.  To
continue the preparation for the proof of
Proposition~\ref{prop:finitesetavoidance} we discuss the pull-back of
the logarithmic tangent sheaf via a general morphism in $\mathcal H'$.

\begin{lem}\label{lem:pbglobgen}
  Under the assumptions of Proposition~\ref{prop:finitesetavoidance}, let $f\in
  \mathcal H'$ be a general element. Then $f^*\bigl(T_S(-\log D)\bigr)$ is globally
  generated on $\P^1$.
\end{lem}
\begin{proof}
  Set $n := \dim S$. Working on $\P^1$, it suffices to show that $f^*(T_S(-\log D))$
  is generated by global sections at a general point $y \in \P^1$, i.e., that there
  exist sections $\sigma_1, \ldots \sigma_n \in H^0(\P^1, f^*(T_S(-\log D)))$ that
  are linearly independent at $y$.
      
  To construct $\sigma_1$, observe that the natural action of $\C^*$ on $\P^1$ that
  fixes $x_0 = [0:1]$ and $x_{\infty} = [1:0]$ yields a non-trivial deformation of
  $f$ in $\mathcal H$. Let $\sigma_1$ be an associated infinitesimal deformation
  which, by general choice of $y$, does not vanish at $y$.
      
  In order to find $\sigma_2, \ldots, \sigma_n$, observe that the curves associated
  with $\mathcal H$ dominate $S$. By general choice of $f$, we can therefore assume
  that the universal morphism
  $$
  \mu: \P^1 \times \mathcal H \to S
  $$
  has rank $n$ at $(y,f)$. The description \cite[II.3.4]{K96} of the tangent
  morphism $T\mu$ then yields the existence of infinitesimal deformations $\sigma_2,
  \ldots, \sigma_n$ whose evaluations $\sigma_i(y)$ along with $\sigma_1(y)$ are
  linearly independent and not tangent to the image of $f$.
\end{proof}
  
\begin{cor}\label{cor:hissmooth}
  Under the conditions of Lemma~\ref{lem:pbglobgen}, both $\Hom_{\bir}(\P^1, S)$ and
  $\mathcal H'$ are reduced and smooth at the point $f$.
\end{cor}
\begin{proof}
  Lemma~\ref{lem:pbglobgen} implies that $f^*(T_S)$ is also globally generated on
  $\P^1$ since it contains the globally generated locally free subsheaf
  $f^*(T_S(-\log D))$ of the same rank. Then $H^1(\P^1, f^*(T_S))=0$, so
  $\Hom_{\bir}(\P^1, S)$ is reduced and smooth of dimension $h^0(\P^1, f^*(T_S))$ by
  \cite[I.2.16]{K96}. This, combined with Lemma~\ref{lem:dimofH} implies that 
  \begin{equation}
    \label{eq:egy}
      h^0(\P^1, f^*(T_S)) - d\leq \dim_f\mathcal H'\leq \dim T_{\mathcal H'}\resto f.
  \end{equation}
  Proposition~\ref{prop:tgtspace} and the fact that $\deg_{\P^1} f^*(T_S) =
  \deg_{\P^1} f^*(T_S(-\log D)) + d$ imply that 
  \begin{equation}
    \label{eq:ket}
  h^0(\P^1, f^*(T_S)) - d=h^0(\P^1, f^*(T_S(-\log D)))= \dim T_{\mathcal H'}\resto f,    
  \end{equation}
  The (in)equalities \eqref{eq:egy} and \eqref{eq:ket} together imply that $\dim T_{\mathcal H'}\resto f
  = \dim_f \mathcal H'.$ Therefore, we obtain that $\mathcal H'$ is reduced and
  smooth at the point $f$.
\end{proof}

\begin{proof}[Proof of Proposition~\ref{prop:finitesetavoidance}]
  \label{proof:finitesetavoidance}
  Consider the standard diagram
  $$
  \xymatrix{
    \P^1 \times \mathcal H' \ar[rr]^{\mu}_{\ \ \text{ univ.~morphism}}
    \ar[d]_{\pi}^{\text{projection}} && S \\ 
    \mathcal H', }
  $$
  and let $\mathcal M := \bigl(\mu^{-1}(M)\bigr)_{\red}$ be the
  set-theoretic preimage of $M$ via $\mu$. Since $\pi$ is proper, it
  is enough to prove that $\pi(\mathcal M)\neq \mathcal H'$.  Assume
  to the contrary, i.e., assume that $\mathcal M$ surjects onto
  $\mathcal H'$ and choose a point $y \in (\pi\resto{\mathcal
    M})^{-1}(f)$.  Since $\mathcal H'$ is smooth at $f$, the general
  choice of $f$ implies that $\pi|_{\mathcal M}$ is \'etale at
  $(y,f)$.  Then the global generation of $f^*\bigl(T_S(-\log
  D)\bigr)$ and the standard description of the tangent morphism
  $T\mu$, \cite[II.3.4]{K96}, yield that the rank of $T\mu|_{\mathcal
    M}$ at $(y,f)$ is at least $n-1$. In particular, $\codim_S M \leq
  1$, a contradiction.
\end{proof}

\subsection{Logarithmic differentials}

Throughout the proof of the main theorem we need to use the sheaf
$\Omega^1_S(\log D)$ of 1-forms with logarithmic poles along $D$. For
the definition and detailed discussion of this notion the reader is
referred to either \cite[Chap.~3]{Deligne70} or \cite[§~2]{EV92}. We
will need to describe $\Omega^1_S(\log D)$ in terms of its restriction
to curves in $S$.

\begin{lem}\label{lem:restrofdiffs}
  Let $F \subset S$ be a smooth curve that intersects $D$
  transversally.  Then the restriction $\Omega^1_S(\log D)\resto F$ is
  an extension of line bundles, as follows:
  \begin{equation}
    \label{eq:restoflogomega1}
    0 \to N^\vee_{F/S} \to \Omega^1_S(\log D)\resto F \to
    \Omega^1_{F}(\log D\resto F) \to 0.
  \end{equation}
  If $D=\sum_{i=1}^rD_i$ is the decomposition of $D$ to irreducible
  components, then the restriction $\Omega^1_S(\log D)\resto {D_1}$ is
  an extension of line bundles, as follows:
  \begin{equation}
    \label{eq:restoflogomega2}
    0 \to \Omega^1_{D_1}(\log (D-D_1)\resto {D_1}) \to \Omega^1_S(\log
    D)\resto {D_1} \to \O_{D_1} \to 0.
  \end{equation}
  Furthermore, if $\dim S = 2$, then
  \begin{equation}
    \label{eq:restoflogomega3}
    \Omega^1_{D_1}(\log (D-D_1)\resto {D_1}) \simeq
      \Omega^1_{D_1} \otimes \O_{D_1}\bigl((D-D_1)\resto {D_1}\bigr).
  \end{equation}
\end{lem}

\begin{proof}
  \PreprintAndPublication{To prove~\eqref{eq:restoflogomega1}, consider the following
    diagram with exact rows \cite[2.3a]{EV92}:
    $$
    \xymatrix{ 0 \ar[r] & \Omega^1_S\resto F \ar[r]
      \ar[d]^{\theta_1} & \Omega^1_S(\log D)\resto F
      \ar[r]\ar[d]^{\theta_2} & \bigoplus_{i=1}^r \O_{D_i}\resto F
      \ar[r]\ar[d]^{\theta_3}  & 0\\
      0 \ar[r] & \Omega^1_F \ar[r] & \Omega^1_F(\log D\resto F) \ar[r]
      & \bigoplus_{i=1}^r \O_{D_i\resto F} \ar[r] & 0. }
    $$
    Notice that $\theta_3$ is an isomorphism and $\theta_1$ is
    surjective, and hence it follows from the Snake Lemma that
    $\theta_2$ is also surjective and $\ker\theta_2\simeq \ker
    \theta_1\simeq N^\vee_{F/S}$. This
    shows~\eqref{eq:restoflogomega1}.
    
    To prove \eqref{eq:restoflogomega2}, consider the following diagram with exact
    rows \cite[2.3c]{EV92}:
    $$
    \xymatrix{ \Omega^1_S(\log D)(-D_1)\ \ar@{^{(}->}[r]
      \ar[d]^{\vartheta_1} & \Omega^1_S(\log (D-D_1))
      \ar@{->>}[r]\ar[d]^{\vartheta_2} &
      \Omega^1_{D_1}(\log(D-D_1)\resto {D_1}) \ar@{.>}[d]^{\vartheta_3}  \\
      \Omega^1_S(\log D)(-D_1)\ \ar@{^{(}->}[r] & \Omega^1_S(\log D)
      \ar@{->>}[r] & \Omega^1_{S}(\log D)\resto {D_1}. }
    $$
    Since the rows are exact, the morphisms $\vartheta_1$ and $\vartheta_2$ imply
    the existence of $\vartheta_3$.  Observe that $\vartheta_1$ is an isomorphism and
    $\vartheta_2$ is injective, hence it follows from the Snake Lemma that
    $\vartheta_3$ is also injective and $\coker\vartheta_3\simeq \coker
    \vartheta_2\simeq \O_{D_1}$ \cite[2.3b]{EV92}.
    
    Finally, if $\dim S = 2$, then $D_1$ is a smooth curve,
    and~\eqref{eq:restoflogomega3} follows immediately from the
    definition.}{\eqref{eq:restoflogomega1} and
    \eqref{eq:restoflogomega2} follow from \cite[2.3a, 2.3c]{EV92}
    using the Snake Lemma.
    
    If $\dim S = 2$, then $D_1$ is a smooth curve,
    and~\eqref{eq:restoflogomega3} follows from the definition.}
\end{proof}

\section{Controlled minimal models of birationally ruled surfaces}
\label{subsec:rel-min}

In this section, we consider log pairs $(S,D)$, where $S$ is a birationally ruled
surface whose boundary intersects the ruling with multiplicity two. More precisely,
we make the following assumption throughout the present section.

\begin{assumption}\label{ass:more}
  Let $S$ be a smooth projective surface and $D\subset S$ a simple normal crossing
  divisor. Assume that there exists a morphism $\pi:S\to C$ whose general fiber is
  isomorphic to $\P^1$.  If $t\in C$ is any point, set $S_t := \pi^{-1}(t)$ and
  assume that $D\cdot S_t = 2$.
\end{assumption}


Our principle aim in this section is to relate the logarithmic Kodaira dimension
$\kappa(S\setminus D)$ with the genus of the base curve $C$ and with the number and
type of fiber components contained in $D$.

The relation in question is formulated in
Propositions~\ref{prop:ruledcanondivisor} and \ref{prop:Dh-isolation}
using a certain sequence of blowing down vertical $(-1)$-curves which
simplifies the self-intersection graph of $D$ and eventually leads to
a $\P^1$-bundle over $C$. The construction of this sequence is
explained in Section~\ref{sec:algorithm} below.

\subsection{Construction of a minimal model. Setup of notation}
\label{sec:algorithm}

To describe the sequence of blowings down we use the following terminology.

\begin{notation}
  A curve $E \subset S$ is called \emph{vertical} if it maps to a
  point in $C$.  Let $D=D^h+D^v$ be the associated decomposition of
  the divisor $D$, where $D^v$ is the sum of the vertical components,
  and $D^h$ the components that surject onto $C$.
\end{notation}

Now consider the sequence of blowings down of vertical $(-1)$-curves, as given by
Algorithm~\ref{constr:1} on page~\pageref{constr:1} below.
The construction obviously depends on choices and is therefore not unique.  While the
results stated in section~\ref{sec:algprop} are independent of the choices made, we
fix a particular set of choices for the remainder of the section and do not pursue
the uniqueness question further.

\setcounter{algocf}{\value{thm}}
\stepcounter{thm}
\begin{algorithm}
  \caption{Construction of a good relative minimal model of $S$}
  \label{constr:1} 
  \dontprintsemicolon
  \BlankLine
  \CommentSty{Step 0: Setup}
  \BlankLine
  $i := 0$, \quad $S_0 := S$, \quad $D^h_0 := D^h$, \quad $D^v_0 := D^v$\;
  \BlankLine
  \CommentSty{Step 1: blow down curves that are disjoint from $D^h$}
  \BlankLine
  \While{there exists a vertical $(-1)$-curve $E_i \subset S_i$, disjoint
    from $D^h_i$}{
    $S_{i+1} \, := $ blow-down of $S_i$ along $E_i$\;
    $D^h_{i+1} := $ cycle-theoretic image of $D^h$ in $S_{i+1}$\;
    $D^v_{i+1} := $ cycle-theoretic image of $D^v$ in $S_{i+1}$\;
    $i \quad\,\,\,\,\, \gets i+1$\;
  }
  $\boxed{k_1 := i}$\;
  \BlankLine
  
  \CommentSty{Step 2: for each reducible fiber $F$ blow down $(-1)$-curves contained
    in $F$, always taking curves in $D^v$ if possible. Stop if $D^h_i$ and $D^v_i$ no
    longer intersect in $F$.}

  \BlankLine
  \For{each reducible fiber $F \subset S_i$}{
    \While{$D^h_i \cap D^v_i \cap F \not = \emptyset$ and there exists a $(-1)$-curve
      in $F$ }{ 
      \eIf{there exists a vertical $(-1)$-curve in $F$, contained in $D^v_i$}{
        $E_i := $ a vertical $(-1)$-curve in $F$, contained in $D^v_i$\;
      }{
        $E_i := $ any vertical $(-1)$-curve in $F$\;
      }
      $S_{i+1} \, := $ blow-down of $S_i$ along  $E_i$\;
      $D^h_{i+1} := $ cycle-theoretic image of $D^h$ in $S_{i+1}$\;
      $D^v_{i+1} := $ cycle-theoretic image of $D^v$ in $S_{i+1}$\;
      $i  \quad\,\,\,\,\, \gets i+1$\;
    }
  }
  $\boxed{k_2 := i}$\;
  \BlankLine
  \CommentSty{Step 3: blow down the remaining vertical $(-1)$-curves}
  \BlankLine
  \While{there exists a vertical $(-1)$-curve $E_i \subset S_i$}{
    $S_{i+1} \, := $ blow-down of $S_i$ along  $E_i$\;
    $i \quad\,\,\,\,\,  \gets i+1$\;
  }
  $\boxed{m := i}$\;
  \BlankLine
  \CommentSty{Now $S_i$ does not contain any vertical $(-1)$-curve, and is therefore
    relatively minimal over $C$.} 
\end{algorithm}

\begin{notation}\label{not:blowdownnotation} \label{diag:blowdown}
  We fix a set of choices, set $S = S_0$ and denote the morphisms that
  occur in Algorithm~\ref{constr:1} as follows.
  $$
  \xymatrix{ S_0 \ar@/^.5cm/[rrrr]^{\rho_i}
    \ar[rr]_{\beta_0:\text{ blow-down}} & & S_1 \cdots
    \ar[rr]_{\beta_{i-1}: \text{ blow-down}} & & S_i
    \ar[rr]_(.45){\beta_i: \text{ blow-down}}
    \ar@/^.5cm/[rrrr]^{\pi_i} & & \cdots S_m
    \ar[rr]^{\pi_m}_{\P^1-\text{bundle}} & & C. }
  $$
  If $t \in C$ is any point, let $S_t := \pi^{-1}(t)$ and $S_{i,t}
  := \pi_i^{-1}(t)$ be the scheme-theoretic fibers. In addition, we
  will also consider the following objects.
  $$
  \begin{array}{lcl}
    k_1, k_2, m & \ldots & \text{the indexes marking the end of Steps~1,~2,~and~3 in
      Algorithm~\ref{constr:1}}\\
    E_i & \ldots & \text{the $\beta_i$-exceptional vertical $(-1)$-curve in $S_i$} \\
    D^h_i, D^v_i & \ldots & \text{the cycle-theoretic images of $D^h$ and
      $D^v$  in $S_i$, respectively} \\
    C_0 & \ldots & \text{the section of $\pi_m$ with minimal self-intersection number} \\
    F_m & \ldots & \text{the numerical class of a fiber of $\pi_m$} \\
    e & \ldots  & - C_0^2, \text{ invariant of the ruled surface $S_m$}\\
    \delta & \ldots & D^h_m \cdot C_0, \text{intersection number of $D^h_m$ and
      $C_0$} \\
  \end{array}
  $$
\end{notation}

\subsection{Properties of the construction}
\label{sec:algprop}

The following two propositions that describe features of the morphisms defined in
\eqref{diag:blowdown} will be shown in Section~\ref{sec:proofofAlgProperties} below.

The first proposition gives a formula for the numerical class of the log canonical
bundle.  This is later used in Section~\ref{sec:K0} to give a relation between the
logarithmic Kodaira dimension of $S\setminus D$, the genus of the base curve and the
number of fibers contained in $D^v$.

\begin{prop}\label{prop:ruledcanondivisor}
  There exists an effective divisor $E' \subset S$, whose support is
  exactly the exceptional locus of $\rho_{k_1}: S \to S_{k_1}$, such
  that the following equality of numerical classes holds.
  $$
  K_S+D\equiv (e+\delta+2g(C)-2)\rho_m^*(F_m)+D^v+E'.
  $$
\end{prop}

\noindent
We will later be interested in reducing to a situation where the horizontal
components are isolated in $D$. The second proposition gives a criterion that
together with Proposition~\ref{prop:ruledcanondivisor} can be used to guarantee that
$D^h_{k_2}$ and $D^v_{k_2}$ intersect only in a controllable manner, if at all.

\begin{prop}\label{prop:Dh-isolation}
  Using the notation of Proposition~\ref{prop:ruledcanondivisor}, let
  $t\in C$ be a point such that the set-theoretic fiber $\supp (S_t)$
  is not contained in the support of $D^v+E'$.  Then $D^h_{k_2}$ and
  $D^v_{k_2}$ do not intersect over $t$, i.e., $t \not \in
  \pi_{k_2}(D^h_{k_2} \cap D^v_{k_2})$.
\end{prop}

\subsection{Proofs of Propositions~3.5 and
  3.6} \label{sec:proofofAlgProperties}

The proofs are not very complicated. They do, however, require some
preliminary computations.

\begin{lem}\label{lem:-1-curves}
  Let $t \in C$ be a point and $i < m$ a number such that $S_{i,t}$ is
  reducible. Then either $S_{i,t}$ contains at least two
  $(-1)$-curves, or it contains exactly one, but with multiplicity
  more than one.
\end{lem}
\begin{proof}
  By blowing down vertical $(-1)$-curves disjoint from $S_{i,t}$, we can assume
  without loss of generality that $i = 0$, and that all vertical $(-1)$-curves blown
  down in Algorithm~\ref{constr:1} lie over $t$. We will then prove the statement by
  induction on $m-i$:
  
  \noindent
  {\it Start of induction, $i=m-1$} In this case, $S_{m-1,t}$ contains exactly two
  $(-1)$ curves.
  
  \noindent
  {\it Induction step} Suppose $i < m-1$ and assume that the statement holds for
  $S_{i+1}$.  Set $x := \beta_i(E_i) \subset S_{i+1}$. Then there are three
  possibilities:
  \begin{enumerate} 
  \item The point $x$ is contained in two vertical $(-1)$-curves. In
    this case, the curve $E_i$ appears in $S_{i,t}$ with multiplicity
    more than one.
    
  \item The point $x$ is contained in exactly one vertical
    $(-1)$-curve $E \subset S_{i+1,t}$. In this case the number of
    $(-1)$-curves in $S_{i,t}$ equals the number of $(-1)$-curves in
    $S_{i+1,t}$, and the multiplicity of $E_i$ in $S_{i,t}$ is at
    least the multiplicity of $E$ in $S_{i+1,t}$.
    
  \item The point $x$ is not contained in a vertical $(-1)$-curve.
    Then $S_{i,t}$ contains at least two $(-1)$ curves.
  \end{enumerate}
  In either case, the claim is shown. This ends the proof of
  Lemma~\ref{lem:-1-curves}.
\end{proof}

\begin{cor}\label{cor:constintersect}
  Let $k_1 \leq j \leq m$ and $E\subseteq S_j$ a vertical $(-1)$-curve.  Then $D^h_j
  \cdot E= 1$.
\end{cor}
\begin{proof}
  Let $\beta := \beta_j\circ\beta_{j-1}\circ\dots\circ\beta_{k_1}: S_{k_1}\to S_j$.
  By construction, $\beta^{-1}(E)$ contains a $(-1)$-curve, $E'\subseteq \beta^{-1}(E)
  \subset S_{k_1}$. By definition of $k_1$, we have $D_{k_1}^h\cap E'\neq
  \emptyset$. Therefore $D_j^h\cap E\supseteq \beta(D_{k_1}^h\cap E')\neq \emptyset$
  and hence
  \begin{equation}
    \label{eq:geqone}
    D^h_j \cdot E\geq 1.
  \end{equation}
  Note that \eqref{eq:geqone} holds for \emph{any} vertical $(-1)$-curve in $S_j$.
  
  Let $S_{j,t}$ be the fiber containing $E$. By Assumption~\ref{ass:more} $D_j^h
  \cdot E \leq D_j^h \cdot S_{j,t} = 2$. Assume that $D_j^h \cdot E = 2$. Then the
  multiplicity of $E$ in $S_{j,t}$ must be one and $D_j^h$ cannot intersect any
  component of $S_{j,t}$ other than $E$. But by Lemma~\ref{lem:-1-curves}, there has
  to be another vertical $(-1)$-curve $E''\subset S_{j,t}$, which is then disjoint
  from $D_j^h$.  This, however, contradicts \eqref{eq:geqone} applied for $E''$,
  hence $D_j^h \cdot E < 2$ and the statement follows.
\end{proof}

\begin{num}[\it Proof of Proposition~\ref{prop:ruledcanondivisor}]\rm 
  The classical formula for the canonical bundle of a blow-up surface
  states:
  $$
  K_{S_i} \equiv \beta_i^* (K_{S_{i+1}}) + E_i,
  $$
  Next we wish to express $D_i^h$ in terms of $E_i$ and the pull-back of
  $D^h_{i+1}$. Depending on $i$, there are two possibilities:
  \begin{description}
  \item[$\boldsymbol{i < k_1}$] In this case $E_i$ and $D^h_i$ are disjoint by
    construction, so $D^h_i \equiv \beta_i^*(D^h_{i+1})$ and therefore
    $$
    K_{S_i}+D^h_i \equiv \beta_i^* (K_{S_{i+1}}+D^h_{i+1}) + E_i.
    $$
    
  \item[$\boldsymbol{i \geq k_1}$] In this case $D^h_i\cdot E_i= 1$ by
    Corollary~\ref{cor:constintersect}, hence $D^h_i \equiv \beta_i^*(D^h_{i+1})-E_i$
    and therefore
    $$
    K_{S_i}+D^h_i \equiv \beta_i^* (K_{S_{i+1}}+D^h_{i+1}).
    $$
  \end{description}
  In summary, we have
  \begin{equation}
    \label{eq:summary}
    K_S+D \equiv \rho^* (K_{S_m}+D^h_m) + D^v + E',    
  \end{equation}
  where $E'$ is an effective divisor supported on the exceptional locus of the
  morphism $\rho_{k_1} : S \to S_{k_1}$.  The standard formula
  \cite[V.~Cor.~2.11]{Ha77} for the canonical bundle of a ruled surface and a simple
  intersection number calculation yields that
  $$
  K_{S_m}\equiv -2 C_0 + (2g(C)-2-e) F_m \text{\quad and \quad } D_m^h \equiv 2
  C_0+(\delta+ 2 e)F_m.
  $$
  Combined with \eqref{eq:summary} this finishes the proof of
  Proposition~\ref{prop:ruledcanondivisor}. \qed
\end{num}

\smallskip

\begin{num}[\it Proof of Proposition~\ref{prop:Dh-isolation}]\rm 
  Let $t \in C$ be a point as in the statement of
  Proposition~\ref{prop:Dh-isolation}.  Assume to the contrary, i.e., that $t \in
  \pi_{k_2}(D^h_{k_2} \cap D^v_{k_2})$. Observe that with this assumption the
  ``while'' condition in Step~2 of Algorithm~\ref{constr:1} stopped only because
  there were no further $(-1)$-curves in the fiber over $t$.  This implies that
  $S_{k_2, t}$ is reduced, irreducible and contained in $D^v_{k_2}$:
  \begin{equation}
    \label{eq:tristan}
    S_{k_2, t} \simeq \P^1  \text{\quad and \quad} S_{k_2, t}
    \subseteq \supp \bigl( (\rho_{k_2})_*(D^v) \bigr) =
    \supp \bigl( (\rho_{k_2})_*(D^v+E') \bigr).
  \end{equation}
  In contrast to~\eqref{eq:tristan}, since $E'$ is supported exactly on the
  exceptional locus of $\rho_{k_1}$, the assumption of
  Proposition~\ref{prop:Dh-isolation} says precisely that
  \begin{equation}
    \label{eq:isolde}
    \supp(S_{k_1, t}) \not \subset \supp \bigl( (\rho_{k_1})_*(D^v)  
  \bigr) = \supp \bigl( (\rho_{k_1})_*(D^v+E')
  \bigr).
  \end{equation}
  
  Now let
  \begin{equation}
    \label{eq:maxdef}
      j := \max \left\{ i \,\,|\,\, \supp(S_{i,t}) \not \subset \supp
    \bigl( (\rho_i)_* (D^v+E') \bigr) \right\}.
  \end{equation}
  It follows by \eqref{eq:tristan} and \eqref{eq:isolde} that
  \begin{equation}\label{eq:rangeofj}
    k_1 \leq j < k_2.    
  \end{equation}
  Loosely speaking, the exceptional curve $E_j$ is the last $(-1)$-curve contracted
  over $t$ that is not in the image of $D^v+E'$. The following two statements follow
  immediately from the choice of $j$.
  \begin{enumerate-c}
  \item \ilabel{castor} $E_j$ is contained in the fiber over $t$, i.e., $E_j \subset
    S_{j,t}$
  \item \ilabel{pollux} $E_j$ is not contained in the image of
    $D^v+E'$, i.e.,
    $E_j \not \subset \supp \bigl( (\rho_j)_* (D^v) \bigr)$.
  \end{enumerate-c}
  The choice of $j$ and the ``if'' statement in Step~2 of Algorithm~\ref{constr:1}
  guarantee that $E_j$ is the only $(-1)$-curve contained in $S_{j,t}$.  In that case
  Lemma~\ref{lem:-1-curves} asserts that the multiplicity of $E_j$ in $S_{j,t}$ is at
  least $2$.  In addition, the first inequality of~\eqref{eq:rangeofj} and
  Corollary~\ref{cor:constintersect} assert that $D^h_j$ intersects $E_j$
  non-trivially. Then by Assumption~\ref{ass:more} $D^h_j$ does not intersect any
  component of the fiber $S_{j,t}$ other than $E_j$.  Then \iref{castor} and
  \iref{pollux} imply that
  $$
  D^h_j \cap D^v_j \cap S_{j,t} = \emptyset.
  $$
  This, combined with the second inequality in~\eqref{eq:rangeofj} above
  contradicts the choice of $k_2$ as the index marking the end of Step~2 of
  Algorithm~\ref{constr:1}.
\end{num}

\part{PROOF OF THE MAIN THEOREM}

\section{Setup of notation}
\label{sec:setup}

In this section, we briefly fix notation used throughout the proof of
Theorem~\ref{thm:main}. The proof will be given in
Sections~\ref{sec:k1}--\ref{sec:Kinfinity} for the cases when the logarithmic Kodaira
dimension of $S^\circ$ is 1, 0 or $-\infty$, respectively. As one might expect, the
case of $\kappa(S^\circ)=0$ is by far the longest and most involved.

\begin{notation}\label{not:main}
  Throughout the rest of the article, we keep the notation and
  assumptions of Theorem~\ref{thm:main}. We fix a smooth projective
  compactification $S$ of $S^\circ$ such that $D = S \setminus S^\circ
  \subset S$ is a simple normal crossing divisor.  Furthermore, let
  $X$ be a smooth projective variety and $f:X\to S$ a morphism such
  that $X\setminus f^{-1}(D)\simeq X^\circ$ and
  $f\resto{X^\circ}=f^\circ$.
  
  Part of the argumentation involves the log minimal model of $(S,D)$.  We will
  therefore adhere to the notation introduced in section~\ref{sec:logMMP}. In
  particular, we use
  $$
  \phi : (S,D) \to (S_{\lambda}, D_{\lambda}) 
  $$
  to denote the birational morphism from $S$ to its logarithmic
  minimal model that is described in Fact~\ref{fact:lmmp}.
\end{notation}

\section{Logarithmic Kodaira dimension $1$}
\label{sec:k1}

If $\kappa(S^\circ)=1$, the statement of Theorem~\ref{thm:main}
follows almost immediately from the logarithmic minimal model program.

\begin{proof}[Proof of Theorem~\ref{thm:main} when $\kappa({S^\circ})=1$]
  By Fact~\ref{fact:lmmp}, we can run the logarithmic minimal model
  program and find a birational morphism $\phi : S \to S_{\lambda}$
  from $S$ to a normal surface $S_{\lambda}$ such that the associated
  log-canonical divisor $K_{S_{\lambda}} + D_{\lambda}$ on
  $S_{\lambda}$ is nef.
    
  The logarithmic abundance theorem in dimension 2,
  Fact~\ref{fact:abundance}, then asserts that for $n \gg 0$ the
  linear system $|n(K_{S_{\lambda}}+D_{\lambda})|$ yields a morphism
  to a curve $\pi_\lambda : S_\lambda \to C$, such that
  $K_{S_{\lambda}} + D_{\lambda}$ is trivial on the general fiber
  $F_\lambda$ of $\pi_\lambda$. Likewise, if $\pi := \pi_\lambda \circ
  \phi$, and $F \subset S$ is a general fiber of $\pi$, then $K_S + D$
  is trivial on $F$.  It follows that $F$ is either an elliptic curve
  that does not intersect $D$, or that $F$ is a rational curve that
  intersects $D$ in two points. It follows, in the former case from
  \cite[Thm.\ 1, Cor.\ 3.2]{Kovacs96e} and in the latter case from
  \cite[0.2]{Kovacs00}, that $f^\circ$ is isotrivial over $F \setminus
  D$, and therefore $\Var(f^\circ) \leq 1 = \kappa(S^\circ)$.
\end{proof}

\section{Logarithmic Kodaira dimension $0$}
\label{sec:K0}

Throughout the present section, we maintain the notation and assumptions of
Theorem~\ref{thm:main} and Section~\ref{sec:setup} and assume that
$\kappa({S^\circ})=0$.

As the proof is rather long, we subdivide it into several steps. We start in
Section~\ref{subsec:ViehwegZuo} by recalling a result of Viehweg and Zuo on which
much of the argumentation is based.  As a first application, we will in
Section~\ref{sec:reduction-to-uniruled} reduce to the situation where $S$ is
uniruled. In Section~\ref{sec:reduction-to-birationally-ruled} we will further reduce
to the case where $S$ is birationally ruled over a curve.

This makes it possible in Section~\ref{subsec:construction-of-good-model} to employ
the results of Chapter~\ref{subsec:rel-min} to construct a birational model of $S$ to
which the aforementioned result of Viehweg and Zuo can be applied. The application
itself, carried out in
Sections~\ref{subsec:2ndapplication}--\ref{subsec:8-computation}, shows that
$\Var(f^\circ) = 0$ and finishes the proof of Theorem~\ref{thm:main}.

\subsection{A result of Viehweg and Zuo}
\label{subsec:ViehwegZuo}

The argumentation relies on the following result describing the sheaf of logarithmic
differentials on the base of a family of canonically polarized varieties.  Note that
we are still using Notation~\ref{not:main}.

\begin{thm}[\protect{\cite[Thm.~1.4(i)]{VZ02}}]\label{thm:VZ}
  There exists an integer $n > 0$ and an invertible subsheaf $\sA \subset \Sym^n
  \Omega^1_S(\log D)$ of Kodaira dimension $\kappa(\sA) \geq \Var(f^\circ)$. \qed
\end{thm}

We will show that $\Var(f^\circ)=0$ by a detailed analysis of $\Omega^1_S(\log D)$.
Essentially, we prove that for all numbers $n$ and locally free subsheaves $\sA
\subset \Sym^n \Omega^1_S(\log D)$, the Kodaira dimension of $\sA$ is never positive,
$\kappa(\sA) \leq 0$.

\subsection{Reduction to the uniruled case} 
\label{sec:reduction-to-uniruled}

A surface $S$ with $\kappa(S^\circ)=0$, of course, need not be uniruled. Using the
result of Viehweg and Zuo, however, we will show that any family of canonically
polarized varieties over a non-uniruled surface $S$ with $\kappa(S^\circ)=0$ is
isotrivial.

\begin{prop}\label{prop:reduction-to-uniruled}
  If $S$ is not uniruled, i.e., if $\kappa(S) \geq 0$, then
  $\Var(f^\circ)=0$.
\end{prop}

We prove Proposition~\ref{prop:reduction-to-uniruled} using two lemmata.

\begin{lem}\label{lem:ntrivksmin}
  If $n \in \mathbb N$ is sufficiently large and divisible, then
  \begin{equation}
    \label{eq:logmin}
    \O_{S_{\lambda}}(n(K_{S_{\lambda}}+D_{\lambda})) = \O_{S_{\lambda}}.  
  \end{equation}
  In particular, the log canonical $\mathbb Q$-divisor $K_{S_{\lambda}}+D_{\lambda}$
  is numerically trivial.
\end{lem}
\begin{proof}
  \eqref{eq:logmin} is an immediate consequence of the assumption $\kappa(S^\circ)=0$
  and the logarithmic abundance theorem in dimension 2, Fact~\ref{fact:abundance},
  which asserts that the linear system $|n(K_{S_{\lambda}}+D_{\lambda})|$ is
  basepoint-free.
\end{proof}

\begin{lem}\label{lem:kntriv}
  If $\kappa(S) \geq 0$, then $S_{\lambda}$ is $\mathbb
  Q$-Gorenstein, $K_{S_{\lambda}}$ is numerically trivial and
  $D_{\lambda} = \emptyset$.
\end{lem}
\begin{proof}
  Lemma~\ref{lem:ntrivksmin} together with the assumption that $|nK_S|\neq\emptyset$
  for large $n$ imply that $\phi$ contracts all irreducible components of $D$, and
  all divisors in any linear system $|nK_S|$, for all $n \in \bN$. Hence the claim
  follows.
\end{proof}

\begin{proof}[Proof of Proposition~\ref{prop:reduction-to-uniruled}]
  We argue by contradiction and assume to the contrary that both $\kappa(S) \geq 0$
  and $\Var(f^\circ) \geq 1$. Let $H \in \Pic(S_{\lambda})$ be an arbitrary ample
  line bundle.

  \begin{narrow}
    \begin{subclaim}\label{claim:unstability}
      The reflexive sheaf of differentials $(\Omega^1_{S_{\lambda}})^{\vee \vee}$ has
      slope $\mu_H\bigl((\Omega^1_{S_{\lambda}})^{\vee \vee}\bigr) = 0$, but it is
      not semistable with respect to $H$.
    \end{subclaim}
    \begin{proof}[Proof of Claim~\ref{claim:unstability}]
      Fix a sufficiently large number $m > 0$ and a general curve $C_{\lambda} \in |m
      H|$.  Flenner's variant of the Mehta-Ramanathan theorem,
      \cite[Thm.~1.2]{Flenner84}, then ensures that if
      $(\Omega^1_{S_{\lambda}})^{\vee \vee}$ is semistable, then so is its
      restriction $(\Omega^1_{S_{\lambda}})^{\vee\vee}|_{C_{\lambda}}$.
  
      By the general choice, $C_{\lambda}$ is contained in the smooth locus of
      $S_{\lambda}$ and stays off the fundamental points of $\phi^{-1}$.  The
      birational morphism $\phi$ will thus be well-defined and isomorphic along $C :=
      \phi^{-1}(C_{\lambda})$.  Lemma~\ref{lem:kntriv} then asserts that
      $$
      \mu_H \bigl((\Omega^1_{S_{\lambda}})^{\vee \vee}\bigr) =
      \frac{K_{S_{\lambda}} \cdot C_{\lambda}}{2m} = 0,
      $$
      which shows the first claim.
  
      Lemma~\ref{lem:kntriv} further implies that $\codim_{S_{\lambda}} \phi(D) \geq
      2$, and so $C$ is disjoint from $D$. The unstability of
      $(\Omega^1_{S_{\lambda}})^{\vee \vee}$ can therefore be checked using the
      identifications
      \begin{nequation}\label{eq:identlog}
        (\Omega^1_{S_{\lambda}})^{\vee \vee}|_{C_{\lambda}} \cong 
        \left.\Omega^1_{S_{\lambda}}\right|_{C_{\lambda}} \cong
        \left.\Omega^1_S\right|_C \cong \Omega^1_S(\log D)|_C.    
      \end{nequation}
      Since symmetric powers of semistable vector bundles over curves
      are again semistable \cite[Cor.~3.2.10]{HL97}, in order to prove
      Claim~\ref{claim:unstability}, it suffices to show that there
      exists a number $n\in \mathbb N$ such that $\Sym^n
      \Omega^1_S(\log D)|_C$ is not semistable. For that, use the
      identifications~\eqref{eq:identlog} to compute
      \begin{align*}
        \deg_C \Sym^n \Omega^1_S(\log D)|_C & = const^+ \cdot \deg_C
        \left.\Omega^1_S\right|_C \\
        & = const^+ \cdot \deg_{C_{\lambda}} (\Omega^1_{S_{\lambda}})^{\vee
          \vee}|_{C_{\lambda}} && \text{\eqref{eq:identlog}}\\
        & = const^+ \cdot (K_{S_{\lambda}}\cdot C_{\lambda}) = 0. &&
        \text{Lemma~\ref{lem:kntriv}}
      \end{align*}
      Hence, to prove unstability it suffices to show that $\Sym^n
      \Omega^1_S(\log D)|_C$ contains a subsheaf of positive degree.
  
      Theorem~\ref{thm:VZ} implies that there exists an integer $n >
      0$ such that $\Sym^n \Omega^1_S(\log D)$ contains an invertible
      subsheaf $\sA$ of Kodaira dimension $\kappa(\sA)
      \geq 1$. But by general choice of $C_\lambda$, this in turn
      implies that $\deg_C(\sA|_C) > 0$, which shows the
      required unstability.  This ends the proof of
      Claim~\ref{claim:unstability}.
    \end{proof} 
  \end{narrow}
  Claim~\ref{claim:unstability} implies that $\Omega^1_{S_{\lambda}}|_{C_{\lambda}}$
  has a subsheaf of positive degree or, equivalently, that it has a quotient of
  negative degree.  On the other hand, Miyaoka's criterion for uniruledness,
  \cite[Cor.~8.6]{Miy85} or \cite{KST05}, asserts that then $S$ is uniruled, leading
  to a contradiction.
\end{proof}

In view of Proposition~\ref{prop:reduction-to-uniruled}, it suffices to prove
Theorem~\ref{thm:main} under the following additional assumption that we
maintain for the rest of Section~\ref{sec:K0}.

\begin{assumption}  \label{ass:k0-uni}
  In addition to the notation and assumptions introduced above we
  further assume that $S$ is uniruled.
\end{assumption}

\subsection{Reduction to birationally ruled surfaces}
\label{sec:reduction-to-birationally-ruled}

We will now show that $S^\circ$ is dominated by curves that are images of $\mathbb
A^1\setminus \{0\}$. We will then, in Proposition~\ref{prop:C*fibers}, conclude that
unless $f^\circ$ is isotrivial, a general point of $S^\circ$ is contained in exactly
one image of $\mathbb A^1\setminus \{0\}$.  This will exhibit $S$ as a birationally
ruled surface.

\begin{prop}\label{prop:uniLK0} 
  The surface $S$ is dominated by a family of rational curves that
  intersect $D$ in two points, but it is not dominated by rational
  curves intersecting $D$ in one point.
\end{prop}

\begin{rem}
  In Proposition~\ref{prop:uniLK0}, the number of intersection points
  is to be understood in the sense of
  Definition~\ref{defn:nptintersect}.
\end{rem}

\begin{proof}[Proof of Proposition~\ref{prop:uniLK0}]
  Recall from \cite[Thm.~1.1]{KMcK} that $S$ is dominated by rational
  curves that intersect $D$ in one point iff $\kappa(S^\circ) =
  -\infty$, which is not the case.
  \begin{narrow}  
    \begin{subclaim}\label{subclaim:domination}  
      The smooth locus $S_{\lambda} \setminus \Sing(S_{\lambda})$ is
      dominated by rational curves intersecting $D_{\lambda}$ in two
      points.
    \end{subclaim}
    \begin{proof}[{Proof of Claim~\ref{subclaim:domination}}]
      We aim to apply \cite[Prop.~1.4(3)]{KMcK}, and so we need 
      that
      \begin{itemize}
      \item the log canonical divisor $K_{S_{\lambda}}+D_{\lambda}$ is
        numerically trivial, and that
      \item the boundary divisor $D_{\lambda}$ is not empty.
      \end{itemize}
      The numerical triviality of $K_{S_{\lambda}}+D_{\lambda}$ was
      shown in Lemma~\ref{lem:ntrivksmin} above.  To show that
      $D_{\lambda} \not = \emptyset$, we argue by contradiction, and
      assume that $D_\lambda = \emptyset$. 
      Set
      $$
      S_\lambda^1 := S_\lambda \setminus \underbrace{\phi(\text{exceptional set
          of $\phi$})}_{\text{finite, contains $\phi(D)$}}.  
      $$
      Then $S_\lambda^1$ is the complement of a finite set and
      $\phi^{-1}|_{S_\lambda^1}$ is a well-defined open immersion.
      Let $f_\lambda:=\phi\circ f$. Then $X|_{f_\lambda^{-1}(S^1_\lambda)} \to
      S^1_\lambda$ is a smooth 
      family of canonically polarized varieties. Consider the
      following diagram:
      $$
      \xymatrix{
        X \ar[d]_{f_\lambda} & & & **[r] \tilde X := X \times_{S_\lambda} \tilde S
        \ar[lll] \ar[d]^{\tilde f} \\ 
        S_\lambda & & \tilde S_\lambda \ar[ll]_{\alpha}^{\text{\ index-one-cover}}
        & \tilde S \ar[l]_{\beta}^{\text{\ log~resolution}} }
      $$
      where $\alpha$ is the index-one-cover described in
      \cite[5.19]{KM98}, and $\beta$ is the minimal desingularization of
      $\tilde S_\lambda$ composed with blow-ups of smooth points so
      that  $\beta^{-1}(\tilde S_\lambda\setminus \alpha^{-1}(S^1_\lambda))$ is a
      divisor with at most simple normal crossings.
      
      By Lemma~\ref{lem:ntrivksmin}, $K_{S_\lambda}$ is torsion.  Since $\alpha$ is
      étale in codimension one this implies that $K_{\tilde S_\lambda}$ is
      trivial. Furthermore, $\tilde S_\lambda$ has only canonical singularities:
      we have already noted in 
\PreprintAndPublication{Remark~\ref{rem:singsinMMP}}{Fact~\ref{fact:lmmp}}  
      that the singularities of 
      $S_\lambda$ are log-terminal, i.e., they have minimal discrepancy $> -1$. Then
      by \cite[Prop.~5.20]{KM98} the minimal discrepancy of the singularities of
      $\tilde S_\lambda$ is  also $> -1$, and as $K_{\tilde S_\lambda}$ is Cartier,
      the discrepancies actually must be integral and hence $\geq 0$,
      cf.~\cite[proof of Cor.~5.21]{KM98}. Consequently,    
      \begin{nequation}
        \label{star}
        K_{\tilde S} =
        \underbrace{\beta^*(K_{\tilde S_\lambda})}_{\cong \mathcal O_{\tilde S}} +
        (\text{effective and $\beta$-exceptional}).        
      \end{nequation}
      This in turn has two further consequences:
      \begin{enumerate}
      \item[i)] $\kappa(K_{\tilde S}) = 0$. In particular, $\tilde S$ is not
        uniruled.
        
      \item[ii)] If we set $\tilde S^1 := (\alpha \circ \beta)^{-1}(S^1_\lambda)$
        then $\tilde X\resto {\tilde f^{-1}(\tilde S_1)}\to \tilde S^1$ is again a
        smooth family of canonically polarized varieties. Letting  $\tilde D :=
        \tilde S \setminus \tilde S^1$ then $\tilde D$ is exactly the
        $\beta$-exceptional set, and \eqref{star} implies that
        $$ 
        \kappa(\tilde S^1) = \kappa (\underbrace{K_{\tilde S} + \tilde D}_{\hskip
          -1cm  \text{effective, $\beta$-exceptional} \hskip -1cm}) = 0.   
        $$ 
        In particular, Proposition~\ref{prop:reduction-to-uniruled} applies to
        $\tilde f: \tilde X \to \tilde S$ and shows that $\tilde S$ is uniruled.  
      \end{enumerate}
      This is a contradiction and thus the proof of
      Claim~\ref{subclaim:domination} is complete.
    \end{proof}
  \end{narrow}
  If $\phi(D) \subset D_{\lambda} \cup \Sing(S_{\lambda})$, i.e., if all connected
  components of $D$ are either mapped to singular points, or to divisors,
  Claim~\ref{subclaim:domination} and Remark~\ref{rem:liftingOfCurves} immediately
  imply Proposition~\ref{prop:uniLK0}. Likewise, if $S_{\lambda}$ were smooth,
  Proposition~\ref{prop:finitesetavoidance} on small set avoidance would imply that
  almost all curves in the family stay off the isolated zero-dimensional components
  of $\phi(D)$, and Proposition~\ref{prop:uniLK0} would again hold. In the general
  case, when $S_{\lambda}$ is singular, and $d_1, \ldots, d_r$ are smooth points of
  $S_{\lambda}$ that appear as connected components of $\phi(D)$, a little more care
  is required.
  
  If $D'$ is the union of connected components of $D$ which are contracted to the set
  of points $\{d_1, \ldots, d_r\} \subset S_{\lambda}$, it is clear that the
  birational morphism $\phi: S \to S_{\lambda}$ factors via the contraction of $D'$,
  i.e. there exists a diagram
  $$
  \xymatrix{ {S} \ar[r]_{\alpha} \ar@/^0.3cm/[rr]^{\phi} & {S'} \ar[r]_{\beta} &
    {S_{\lambda}}}
  $$
  where $S'$ is smooth, and $\alpha$ maps the connected components of $D'$ to
  points $d'_1, \ldots, d'_r \in S'$ and is isomorphic outside of $D'$.
  
  Now, if $D'' := D \setminus D'$, the above argument shows that $S'$ is dominated by
  rational curves that intersect $\alpha(D'')$ in two points. Since $S'$ is smooth,
  Proposition~\ref{prop:finitesetavoidance} on small set avoidance applies and shows
  that almost all of these curves do not contain any of the $d'_i$. Therefore, we
  have seen that most of the curves in question intersect $\alpha(D)$ in two points.
  Remark~\ref{rem:liftingOfCurves} then completes the proof.
\end{proof}

\begin{prop}\label{prop:C*fibers}
  Either $\Var(f^\circ)=0$, or there exists a smooth curve $C$ and morphisms
  $$
  \xymatrix{ C & & {\tilde S} \ar[rr]^{\psi}_{\text{birational}}
    \ar[ll]_{\pi}^{\text{birat.~ruling}} & & S }
  $$
  such that 
  \begin{enumerate}
  \item \ilabel{i52:1} $\tilde S$ is a smooth surface and $\tilde D :=
    \psi^{-1}(D)$ is a divisor with simple normal crossing support.
    
  \item \ilabel{i52:3}\label{i52:3}\label{intersect-two} If $t \in C$
    is a general point, then the fiber $\tilde S_t := \pi^{-1}(t)$ is
    isomorphic to $\P^1$, and intersects $\tilde D$ in exactly two
    points. In particular, $\tilde D\cdot \tilde S_t=2$ for all $t\in
    C$.
  \item \ilabel{i52:4} The restriction of $f^\circ$ to any fiber of $\pi \circ
    \psi^{-1}|_{S^\circ}$ is isotrivial.
    
  \item The morphism $\psi$ is birational, and isomorphic over $S^\circ$. In
    particular, $\pi$ induces a fibration $\pi \circ \psi^{-1}|_{S^\circ}:S^\circ \to
    C$.
  \end{enumerate}
\end{prop}

\begin{proof}
  If $\Var(f^\circ) = 0$, there is nothing to prove, so we may assume that
  $\Var(f^\circ) > 0$.
  
  By Proposition~\ref{prop:uniLK0}, there exists a proper curve $C'
  \subset \Chow(S)$ such that general points $t \in C'$ are associated
  with irreducible, reduced rational curves $\ell_t$ intersecting $D$
  in exactly two points, in the sense of
  Definition~\ref{defn:nptintersect}. Then by
  \cite[Thm.~0.2]{Kovacs00} the restriction of the family $f$ to a
  general curve $(\ell_t)_{t \in C'}$ is isotrivial, so \iref{i52:4}
  follows from the rest of the statement.
  
  Let $\tilde S'$ be the restriction of the universal family over $\Chow(S)$ to $C'$
  and $\psi': \tilde S' \to S$ the restriction of the cycle morphism.  Finally, let
  $C$ be the normalization of $C'$, $\tilde S$ the normalization of $\tilde
  S'\times_{C'}C$ and $\psi: \tilde S \to S$ the morphism induced by $\psi'$.  After
  blowing up further, we may assume that $\tilde S$ is smooth and that $\tilde D :=
  (\psi^{-1}(D))_{\red}$ has only simple normal crossings.
  
  \iref{i52:1} and \iref{i52:3} hold by construction.  The last part of \iref{i52:3}
  follows from the fact that for a general $t\in C$, $\tilde D$ intersects $\tilde
  S_t$ transversally and the numerical class of $\tilde S_t$ is independent of $t$.
  By Zariski's main theorem \cite[Thm.~V.5.2]{Ha77}, the proof of
  Proposition~\ref{prop:C*fibers} is finished if we show that $\psi$ is birational,
  and that it is finite over $S^\circ = S\setminus D$.
  
  \paragraph{Birationality}
  
  Since we are working in characteristic $0$, it suffices to show that
  $\psi$ is generically injective.  Assume to the contrary that a
  general point in $S^\circ$ is contained in more than one of the
  $\ell_t$'s.
  
  Fix a general $t\in C'$. Then the associated curve $\ell_t$
  intersects $D$ in exactly two points. Further, there exists an open
  set $\ell_t^\circ \subset \ell_t$ such that any $x \in \ell_t^\circ$
  satisfies the following:
  \begin{itemize}
  \item $x$ is a general point of $S^\circ$, and
  \item there exists a point $t_x \in C'$ such that the associated curve $\ell_{t_x}$
    contains $x$, is different from $\ell_t$, and intersects $D$ in exactly two
    points,
  \end{itemize}
  Since the $\ell_{t_x}$ dominate $S$, and since $f$ is isotrivial
  both over $\ell_t$ and over any of the $\ell_{t_x}$, $f$ must be
  isotrivial, contrary to our assumption.
  
  \paragraph{Finiteness}
  
  If there was a point $s \in S^\circ$ that was contained in
  infinitely many of the curves $(\ell_t)_{t \in C'}$, then the
  isotriviality of the restrictions $f|_{\ell_t}$ would again imply
  that $f$ is isotrivial over $S$, contradicting our assumptions.
\end{proof}

To prove Theorem~\ref{thm:main}, we may replace $S$ by $\tilde S$ and $X$ by a
desingularization of $X \times_S \tilde S$.  We will thus make the following
additional assumption that we will maintain without loss of generality for the rest
of the present section.

\begin{assumption}\label{ass:double-section}
  Assume that there exists a morphism $\pi : S \to C$ to a smooth curve $C$, with the
  following property: If $t \in C$ is a general point, then the fiber $S_t :=
  \pi^{-1}(t)$ is isomorphic to $\P^1$, and intersects $D$ in exactly two points. In
  particular, $D\cdot S_t=2$ for all $t\in C$.
\end{assumption}

\subsection{Construction of a good model of $S$}
\label{subsec:construction-of-good-model}

In order to apply Theorem~\ref{thm:VZ} to our setup, we need to study the restriction
of the sheaf of logarithmic differentials to components of the boundary.  While the
restriction to isolated components is easily described using
Lemma~\ref{lem:restrofdiffs}, in general we have very little control over the
intersection graph of the boundary divisor. It seems therefore rather difficult to
describe logarithmic differentials directly in this na\"ive manner.

To overcome this difficulty and to simplify the intersection graph, we recall the
results of Section~\ref{subsec:rel-min} and apply Algorithm~\ref{constr:1} to the
birationally ruled surface $\pi : S \to C$. We have seen in
Proposition~\ref{prop:Dh-isolation} that in the intermediate surface $S_{k_2}$, the
horizontal components of the boundary divisor $D^h_{k_2}$ are disjoint from the
vertical components $D^v_{k_2}$ as long as $D^v+E'$ does not contain an entire fiber
of $\pi$.  This makes the analysis of the sheaf of logarithmic differentials much
easier.  Unfortunately, the image $\rho_{k_2}(D)$ need not be a normal crossing
divisor. We construct a log resolution of $\bigl(S_{k_2}, \rho_{k_2}(D) \bigr)$ as
follows.

\begin{construction}\label{constr:Xmu}
  Let $S^\circ_{k_2} \subset S^\circ$ be the maximal open subset where
  $\rho_{k_2}$ is isomorphic. The difference $S^\circ \setminus
  S^\circ_{k_2}$ is then contained in finitely many fibers.  By
  definition, we can view $S^\circ_{k_2}$ also as an open subset of
  $S_{k_2}$, and observe that
  $$
  S_{k_2} \setminus S^\circ_{k_2} = D_{k_2}^h \cup D_{k_2}^v \cup
  \text{\{finitely many isolated points\}}.
  $$
  Let $\beta : S_{\mu} \to S_{k_2}$ be the minimal log resolution
  of the pair $(S_{k_2}, S_{k_2} \setminus S_{k_2}^\circ)$. If
  $X_{\mu}$ is a desingularization of the pull-back $X
  \times_{S_{k_2}} S_{\mu}$, we obtain a diagram as follows:
  $$
  \begin{array}{ll}
    \xymatrix{
      X_{\mu} \ar[d] \ar[rr]^{f_{\mu}} & & S_{\mu} \ar[d]_{\beta}
      \ar@/^0.3cm/[rrd]^{\pi_{\mu}=\pi_{k_2}\circ\beta} \\
      X \ar[r]_{f} & S \ar[r]_{\rho_{k_2}} & S_{k_2} \ar[r]_{\gamma} & S_m \ar[r]
      _{\pi_m} & C
    }
  \end{array}
  \text{\quad with \quad}
  \begin{array}{ll}
    \pi & = \pi_m \circ \gamma \circ \rho_{k_2} \\
    \pi_{k_2} & = \pi_m \circ \gamma \\
    \pi_{\mu} & =  \pi_m \circ \gamma \circ \beta
  \end{array}
  $$
  Again, the rational map $\beta^{-1}\circ \rho_{k_2}$ is an
  isomorphism over $S^\circ_{k_2}$, so that we can view
  $S^\circ_{k_2}$ as a subset of $S_\mu$. The morphism $f_{\mu}$ is
  smooth over $S^\circ_{k_2} \subset S_\mu$, and to show that
  $f^\circ$ is isotrivial, it suffices to prove isotriviality for
  $f_\mu$.  Finally, let $D_{\mu} := S_{\mu} \setminus S^\circ_{k_2}$.
  Then $D_{\mu}$ is a simple normal crossing divisor that we decompose
  into horizontal and vertical components, $D_{\mu} = D_{\mu}^h \cup
  D_{\mu}^v$, as before. Note that by
  Assumption~\ref{ass:double-section} $D_\mu^h$ is a double section,
  in particular, it has at most two irreducible components.
\end{construction}

\begin{notation}
  We have applied Algorithm~\ref{constr:1} to the birationally ruled
  surface $\pi:S \to C$ in order to construct $S_{k_2}$ and $S_m$.
  Throughout the remainder of the present Section~\ref{sec:K0}, we
  maintain Notation~\ref{not:blowdownnotation} that was introduced on
  page~\pageref{not:blowdownnotation} along with
  Algorithm~\ref{constr:1}.
  
  In particular, we let $C_0\subset S_m$ be the distinguished section of $\pi_m$ with
  the minimal self-intersection number, $e=-C_0^2$ and $\delta=D^h_m\cdot C_0$.
\end{notation}

\subsection{Another application of Theorem~\ref{thm:VZ}}
\label{subsec:2ndapplication}

Fix an irreducible component $D_{\mu}^{h,1} \subset D_{\mu}^h$.  Using
Proposition~\ref{prop:Dh-isolation} we will be able to show in
Section~\ref{subsec:8-computation} below that $D_{\mu}^{h,1}$ is either rational or
elliptic, and compute an upper bound for the number of intersection points between
$D_{\mu}^{h,1}$ and other components of $D_{\mu}$. Theorem~\ref{thm:VZ} will then
apply to $S_{\mu}$ and yield the following proposition.

\begin{prop}\label{prop:restroflog-VZ}
  Let $D_{\mu}^{h,1} \subset D_{\mu}^h$ be an irreducible component. If either one of
  the following holds:
  \begin{enumerate}
  \item $D_{\mu}^{h,1}$ is elliptic and isolated in $D_{\mu}$, or
  \item $D_{\mu}^{h,1}$ is rational and intersects other
    components of $D_{\mu}$ in at most two points,
  \end{enumerate}
  then $\Var(f^\circ) = 0$.
\end{prop}
\begin{proof}
  We argue by contradiction and assume $\Var(f^\circ) > 0$. By
  Theorem~\ref{thm:VZ} there exists a number $n > 0$ and an invertible
  subsheaf $\sA_{\mu} \subset \Sym^n \Omega^1_{S_{\mu}}(\log D_{\mu})$
  of Kodaira dimension $\kappa(\sA_{\mu}) > 0$.
  
  If $F_{\mu} \subset S_{\mu}$ is a general fiber of $\pi_{\mu}$, then $F_{\mu}$ is
  isomorphic to $\P^1$ and intersects $D_{\mu}$ transversally in exactly two points.
  Then the logarithmic normal bundle sequence~\eqref{eq:restoflogomega1} from
  Lemma~\ref{lem:restrofdiffs} is split.  The restriction $\Omega^1_{S_{\mu}}(\log
  D_{\mu})|_{F_{\mu}}$ is therefore trivial, and so is $\Sym^n
  \Omega^1_{S_{\mu}}(\log D_{\mu})|_{F_{\mu}}$. It follows that the restriction of
  $\sA_\mu$ to $F_{\mu}$ is a trivial subsheaf of
  $\Sym^n \Omega^1_{S_{\mu}}(\log D_{\mu})|_{F_\mu}$. This has two consequences.
  First, the restriction of $\sA_\mu$ to $D^{h,1}_\mu$ must have positive Kodaira
  dimension. Second, the natural map between restrictions, 
  $\sA_\mu|_{D^{h,1}_\mu} \to \Sym^n \Omega^1_{S_{\mu}}(\log D_{\mu})|_{D^{h,1}_\mu}$,
  is not zero.
  
  On the other hand, sequence~\eqref{eq:restoflogomega2} from
  Lemma~\ref{lem:restrofdiffs} gives
  $$
  0 \to \underbrace{\Omega^1_{D^{h,1}_{\mu}}(\log (D_{\mu} -
    D^{h,1}_{\mu})|_{D^{h,1}_{\mu}})}_{=: \sL} \to
  \underbrace{\Omega^1_{S_{\mu}}(\log D_{\mu})|_{D^{h,1}_{\mu}}}_{=:
    \sR} \to \O_{D^{h,1}_{\mu}} \to 0,
  $$
  where $\sL \in \Pic(D^{h,1}_{\mu})$ is a line bundle of
  degree
  $$
  \deg \sL = 2g(D^{h,1}_{\mu}) - 2 + \# \{
  \text{intersection points of $D^{h,1}_{\mu}$ with other
    components of $D_{\mu}$}\}.
  $$
  If $D_{\mu}^{h,1}$ is elliptic and isolated in $D_{\mu}$, then $\deg \sL = 0$.
  Then $\sR$ is semistable of degree 0, and so is $\Sym^n \sR$.  Likewise, if
  $D_{\mu}^{h,1}$ is rational and intersects other components of $D_{\mu}$ in at most
  two points, then $-2 \leq \deg \sL \leq 0$. Then $\sR$ is a sum of line bundles of
  non-positive degree, and so is $\Sym^n \sR$.  In both cases, we have
  $\deg(\sA_{\mu}|_{D_{\mu}^{h,1}}) \leq 0$. A contradiction.
\end{proof}

\subsection{Computation of genera and intersection points}
\label{subsec:8-computation}

In order to apply Proposition~\ref{prop:restroflog-VZ}, we need to
compute the genus of $D_{\mu}^{h,1}$ and the number of intersection
points between $D_{\mu}^{h,1}$ and other components of $D_{\mu}$.
While it is possible to write down a (complicated) formula that
involves both pieces of information, we found it easier to consider
the cases where $D_{\mu}^h$ is reducible, respectively irreducible,
separately in Sections~\ref{sec:4red} and \ref{sec:4irr}.

The following simple observation helps to count the number of
intersection points in either case.

\begin{obs}\label{obs:counting-pts}
  If $x \in D_{\mu}^{h,1}$ is a point of intersection between $D_{\mu}^{h,1}$ and
  other components of $D_{\mu}$, then, using the notation introduced in
  \eqref{constr:Xmu}, one of the following holds:
  \begin{enumerate}
  \item \ilabel{ii:1} $\beta(x)$ is a singular point of
    $D_{k_2}^{h,1}=\beta_*(D^{h,1}_\mu)$. In particular, $(\gamma \circ
    \beta)(x)$ is a singular point of
    $D_m^{h,1}=(\gamma\circ\beta)_*(D^{h,1}_\mu)$.
    
  \item \ilabel{ii:3} $\beta(x) \in D^h_{k_2} \cap D^v_{k_2}$. By
    Proposition~\ref{prop:Dh-isolation}, $\pi_{\mu}(x)$ is a point
    whose set-theoretic fiber $\supp(S_{\pi_{\mu}(x)})$ is contained
    in the support of $D^v + E'$. 
    
  \item \ilabel{ii:2} $D^h_{k_2}$ is reducible and $\beta(x)\in D_{k_2}^{h,1} \cap
    D_{k_2}^{h,2}$. In particular, $D^h_m$ is reducible and $(\gamma\circ\beta)(x)\in
    D_m^{h,1} \cap D_m^{h,2}$.  \qed
  \end{enumerate}
\end{obs}

Formulated in more technical terms, Observation~\ref{obs:counting-pts}
gives the following.

\begin{cor}\label{cor:estimate-for-I}
  If we denote the the number of points as follows,
  \begin{align*}
    I & := \# \{\text{intersection points between $D_{\mu}^{h,1}$ and
      other components of
      $D_{\mu}$}\} \\
    I_1 & := \# \{ x \in D_{\mu}^{h,1} \,|\, \text{$(\gamma\circ\beta)(x)$ is a
      singular point of $D_m^{h,1}$}\} \\
    I_2 & := \# \{ x \in D_{\mu}^{h,1} \,|\, \supp(S_{\pi_{\mu}(x)}) \subset
    \supp(D^v + E') \}
  \end{align*}
  then
  $$
  I \leq I_1 + I_2 + \left\{
    \begin{array}{ll}
      D_m^{h,1} \cdot D_m^{h,2} & \text{if $D^h_m$ is reducible}  \\
      0 & \text{otherwise}
    \end{array}
  \right.
  $$
  \qed
\end{cor}

It remains to compute the numbers $I_1$, $I_2$ and $D_m^{h,1} \cdot
D_m^{h,2}$ in all relevant cases. Before we do that, we remark that
the results of Section~\ref{subsec:rel-min} immediately give an upper
bound for the number $I_2$. For this, we maintain the notation of
Section~\ref{subsec:rel-min}. In particular, we use the numbers $e$
and $\delta$ that were introduced in
Notation~\ref{not:blowdownnotation}.

\begin{lem}\label{lem:pwd}
  If $d$ is the degree of the finite morphism $D_{\mu}^{h,1} \to C$, then
  $$
  0 \leq I_2 \leq -d\cdot(e+\delta+2g(C)-2).
  $$
\end{lem}
\begin{proof}
  It is clear from the definition that
  $$
  I_2 \leq d \cdot \# \{ t \in C \,|\, \supp(S_t) \subset \supp(D^v
  + E') \}
  $$
  On the other hand, since $\kappa(S^\circ)=0$, it follows from
  Proposition~\ref{prop:ruledcanondivisor} that
  $$
  0 \leq \# \{ t \in C \,|\, \supp(S_t) \subset \supp(D^v + E') \}
  \leq -(e+\delta+2g(C)-2).
  $$
  This shows the claim.
\end{proof}

\subsubsection{Computation of genera and intersection points if $D_{\mu}^h$ is
  reducible}\label{sec:4red}

Let $D_m^{h,1}$ and $D_m^{h,2}$ be the irreducible components of
$D_m^h$ and write $D_m^{h,i}\equiv C_0+b_iF_m$. Since $\delta = D^h_m\cdot C_0$, a
simple computation shows that $D_m^h \equiv 2C_0 + (\delta+2e)F_m$. In particular,
$b_1+b_2 = \delta +2e$, and then the intersection number between the components is
\begin{newequation}\label{eq:wurtz}
  D_m^{h,1} \cdot D_m^{h,2} = (C_0+b_1F_m)\cdot (C_0+b_2F_m)= -e+b_1+b_2=
  e+\delta. 
\end{newequation}
The number $I$ can be then bounded as follows:
\begin{align*}
  0 \leq I & \leq I_1 + I_2 + D_m^{h,1} \cdot D_m^{h,2} &&
  \text{Corollary~\ref{cor:estimate-for-I}} \\
  &  = I_2 + D_m^{h,1} \cdot D_m^{h,2}  && \text{because $D_m^{h,1} \cong C$ is smooth}\\
  & \leq -(e+\delta+2g(C)-2)+ D_m^{h,1} \cdot D_m^{h,2}  && \text{Lemma~\ref{lem:pwd}}\\
  & = 2-2g(C). && \text{Equation~\eqref{eq:wurtz}}
\end{align*}
In particular, either $D_{\mu}^{h,1}$ is elliptic and $I=0$, or it is
rational and $I\leq 2$.  The prerequisites of
Proposition~\ref{prop:restroflog-VZ} are thus fulfilled in any case.
It follows that $\Var(f^\circ) = 0$, and
Theorem~\ref{thm:main} is shown in this case. \qed

\subsubsection{Computation of genera and intersection points if $D_{\mu}^h$ is irreducible}\label{sec:4irr}

Since $D_m^h \equiv 2C_0 + (\delta+2e)F_m$, the standard formula
\cite[V.~Cor.~2.11]{Ha77} for the numerical class of the canonical
bundle of a ruled surface gives
$$
K_{S_m} + D_m^h \equiv (e+\delta + 2g(C)-2) F_m.
$$
The formula \cite[V.~Ex.~1.3a]{Ha77} for the arithmetic genus of
$D_m^h$ then says that
\begin{newequation}\label{eq:arith_genus}
  p_a (D_m^h) = \frac{(K_{S_m}+D_m^h)\cdot D_m^h}{2}+1 = 
  \underbrace{(e+\delta + 2g(C)-2)}_{\leq 0 
    \text{ by Lemma~\ref{lem:pwd}}} +1 \leq 1    
\end{newequation}
In particular, $D^h_{\mu}$ is either elliptic or rational. We treat
these cases separately.

\bigskip

If $D_{\mu}^h$ is elliptic, then $g(D_m^h) = p_a(D_m^h) = 1$ and so
$D_m^h$ is smooth, $I_1=0$. Corollary~\ref{cor:estimate-for-I} and
Equation~\eqref{eq:arith_genus} assert
$$
I \leq I_2 \leq -2\cdot (e+\delta+2g(C)-2) = -2 \cdot \left(p_a
  (D_m^h)-1\right) = 0.
$$
The elliptic curve $D^h_{\mu}$ is thus isolated in $D_{\mu}$,
Proposition~\ref{prop:restroflog-VZ} implies that $\Var(f^\circ) = 0$,
and Theorem~\ref{thm:main} is shown in this case.

Therefore we may assume that $D_{\mu}^h$ is rational.  If $D^h_m$ is singular, its
arithmetic genus is exactly one and \cite[IV.~Ex.~1.8a]{Ha77} asserts that there are
at most two points in $D^h_{\mu}$ that map to the singularities. Hence, whether
$D^h_m$ is singular or not, $I_1$ is always bounded as $I_1 \leq 2\cdot p_a(D^h_m)$.
Then Corollary~\ref{cor:estimate-for-I}, Lemma~\ref{lem:pwd} and
Equation~\eqref{eq:arith_genus} assert that
$$
I \leq I_1 + I_2  \leq 2\cdot p_a(D^h_m) -2\cdot (e+\delta+2g(C)-2) = 2.
$$
Again, Proposition~\ref{prop:restroflog-VZ} applies and
Theorem~\ref{thm:main} is shown. \qed

\section{Logarithmic Kodaira dimension $-\infty$}
\label{sec:Kinfinity}

We maintain the notation and assumptions of Theorem~\ref{thm:main} and
Section~\ref{sec:setup} and assume that $\kappa({S^\circ})=-\infty$. In this case,
the statement follows quickly from the logarithmic abundance result of Keel-McKernan:

\begin{thm}[{\cite[Thm.~1.1]{KMcK}}]\label{prop:uniLKinfty} 
  Let $S$ be a smooth projective surface and $D \subset S$ a reduced divisor with
  simple normal crossings. Assume that
  $$
  \kappa(S\setminus D) = - \infty.
  $$
  Then $S \setminus D$ is dominated by a family of curves that are isomorphic to
  $\mathbb A^1$. \qed
\end{thm}

\begin{proof}[Proof of Theorem~\ref{thm:main} if $\kappa({S^\circ})=-\infty$] 
  Follows immediately from Theorem~\ref{prop:uniLKinfty} because families over
  $\mathbb A^1$ are necessarily isotrivial by \cite[0.2]{Kovacs00}.
\end{proof}

\section{Generalizations}\label{sec:complements}

\begin{num}\rm
It may be worth to note that the proof of Theorem~\ref{thm:main} uses
only the following two facts about families of canonically polarized
varieties.
\begin{enumerate}
\item\ilabel{il:vz} positive variation guarantees the existence of a
  sheaf $\sA \subset \Sym^n \Omega^1_S(\log D)$ of positive
  Kodaira-Iitaka dimension---see Theorem~\ref{thm:VZ}
\item\ilabel{il:kov} families over $\mathbb P^1$, $\mathbb A^1$, $\mathbb
  A^1\setminus \{0\}$ and over elliptic curves are necessarily trivial.
\end{enumerate}
Since~\iref{il:kov} is an immediate consequence of~\iref{il:vz}, the
proof of Theorem~\ref{thm:main} will work with few modifications
whenever we have a family that guarantees the existence of a subsheaf
of $\Sym^n \Omega^1_S(\log D)$ similar to what we have
in~\iref{il:vz}. For instance, \cite[Thm.~1.4(iii)]{VZ02} applies to
give the following complementary statement to Theorem~\ref{thm:main}.  
\end{num}

\begin{thm}
  Let $S^\circ$ be a smooth quasi-projective surface and $f^\circ:
  X^\circ \to S^\circ$ a smooth family of maximal variation,
  $\Var(f^\circ)=2$ such that $\omega_{X^\circ/S^\circ}$ is relatively
  semi-ample. Then $S^\circ$ is of log general type,
  i.e.,~$\kappa(S^\circ)=2$. \qed
\end{thm}

\begin{rem}
  It is conjectured that $\omega_{X^\circ/S^\circ}$ is relatively semi-ample iff
  $f^\circ$ is a family of minimal manifolds. So far, this is known only if the
  fibers have dimension at most three.
\end{rem}

\begin{cor}
  Viehweg's Conjecture~\ref{conj:viehweg} holds for families of
  minimal curves, surfaces or threefolds over surfaces. \qed
\end{cor}

\providecommand{\bysame}{\leavevmode\hbox to3em{\hrulefill}\thinspace}
\providecommand{\MR}{\relax\ifhmode\unskip\space\fi MR}
\providecommand{\MRhref}[2]{%
  \href{http://www.ams.org/mathscinet-getitem?mr=#1}{#2}
}
\providecommand{\href}[2]{#2}

\end{document}